# A theoretical framework and some promising findings of grey wolf optimizer, part I: analytical model of sampling distribution and stability analysis

Haoxin Wang, Libao Shi*

National Key Laboratory of Power Systems, Shenzhen International Graduate School, Tsinghua University, 518055, P. R. China

**Abstract**

This paper proposes a theoretical framework of the grey wolf optimizer (GWO) based on several interesting theoretical findings, involving sampling distribution, order-1 and order-2 stability, and global convergence analysis. In the part I of the paper, the characteristics of the sampling distribution of the new solution and the probabilistic stability of the GWO are carefully discussed based on the well-known stagnation assumption for simplification purposes. Firstly, the characteristics of the sampling distribution of the new solution, mainly related to the shape of the joint probability density function (PDF), are discussed under the assumption that the original solution before updating is constant. Then, the assumption that the original solution is constant is eliminated to perform the sampling distribution analysis, based on which several characteristics of the new solution are provided, containing the shape of the joint PDF and central moments of any positive integer order. Finally, the order-1 and order-2 stability of the GWO under stagnation assumption is introduced and proved as the inference of conclusions above, which are all verified by numerical simulations.

## 1. Introduction

The last few decades have witnessed the flourishing development of a large class of optimization algorithms called meta-heuristic algorithms, which have been regarded as effective solutions to the increasingly complex optimization problems in real life [1-2]. Different from traditional optimization methods like Conjugate Gradient method or Newton method [3], some random searching patterns and mechanisms were introduced in the meta-heuristic algorithms, which substantially reduced the dependence of algorithm performance on the initial solution and increased the probability of finding the global optimal solution of the optimization problem. Moreover, these meta-heuristic algorithms were proposed mainly based on metaphors of some natural or man-made processes, which typified the idea of bionics [4]. Based on the type of natural or man-made process imitated by the meta-heuristic algorithms, they can usually be classified into the following four categories [5]: evolutionary algorithms inspired from Darwin's theory of evolution [6-8], physics-based algorithms inspired from Laws of physics [9-11], human-based algorithms imitating human behavior [12-14], and swarm-based algorithms imitating the behavior of biological populations [15-18]. For more detailed reviews of the meta-heuristic algorithms, please refer to [2, 19-21].

As a typical swarm-based meta-heuristic algorithm, the grey wolf optimizer (GWO) proposed by S. Mirjalili took its metaphor from the hunting process of grey wolf groups [22]. In nature, the grey wolves mostly prefer to live in groups, in which a grey wolf belongs to one of the four categories: α wolf, β wolf, δ wolf, and ω wolf [23]. During predation, the ω wolves are guided by the α wolf, β wolf, and δ wolf to move towards the prey. In the GWO, each grey wolf in the group was regarded as a search agent, and the α wolf, β wolf, and δ wolf were regarded as the best search agent, the second best search agent, and the third best search agent, respectively. Other search agents (representing ω wolves) updated their positions under the guidance of the positions of α wolf, β wolf, and δ wolf (the detailed updating equation of each search agent can be found in section 2). Since its proposal in 2014, owing to its intuitive procedure and outstanding performance, the GWO has attracted widespread attention and has been applied to solve various optimization problems in real life, such as optimal operation of power systems [24-26], training of neural networks [27-28], clustering analysis [29-30], and path planning [31-32], etc. Please refer to [33-35] for more detailed reviews of applications of the GWO.

Despite the research on applications, studies have been done on the performance analysis of the GWO as well. In fact, the updating mechanism and the distribution of particles within the feasible region were already analyzed qualitatively and presented in the form of diagrams in the original literature of the GWO (please see Fig. 3 – Fig. 5 given in [22]). Later on, the convergence behavior of the GWO was verified and studied by numerical simulations based on various test functions [36], and the performance of the GWO was compared with several well-known meta-heuristic algorithms such as the particle swarm optimization (PSO), differential evolution (DE), firefly algorithm (FA), etc. [37-39]. In addition, the mechanisms that how parameters in the updating equation of the GWO controlled the distribution of particles were studied qualitatively, based on which extensive research has been carried out on the improvement of the original GWO from various aspects, such as parameter selection and control [40-41], revision of the updating equation [42-44], and study of population topology [45], etc. It can be seen that till now, most studies on the performance analysis of the GWO were conducted qualitatively via diagrams or numerical simulations, and there were little research on the theoretical foundations of the GWO based on quantitative analysis, such as sampling distribution, probabilistic stability analysis, local/global convergence analysis, etc.

In the part I of the paper, we propose a theoretical framework of the GWO based on several new theoretical findings, involving sampling distribution, order-1 and order-2 stability, and global convergence analysis. In particular, we analyze the characteristics of the sampling distribution of the new solution, and some promising findings of the GWO are carefully discussed, in which the stagnation assumption is introduced and concerned for simplification. Firstly, the original solution before updating is treated as a constant rather than a random variable. Under this assumption, the new solution could be treated as a function of random parameters. Then, the assumption above is eliminated, that is, the original solution is still a random variable. In this case, the solution should be regarded as a stochastic process corresponding to the iteration of the GWO rather than a simple random variable. Finally, as an application of the sampling distribution characteristics of the solution analyzed above, we introduce and prove the order-1 and order-2 stability of the GWO under stagnation assumption. In the part II of the paper, based on the theoretical results concluded in the part I, we will prove the global convergence of the GWO. All of the conclusions above are verified by numerical simulation results.

The part I of the paper is organized as follows. In section 2, the updating mechanism and procedure of the GWO are introduced in brief. In section 3, the stagnation assumption is introduced, and the sampling distribution of the new solution with both stagnation assumption and constant original solution assumption is solved and analyzed. In section 4, the constant original solution assumption is eliminated, the characteristics of the sampling distribution of the new solution is discussed only under the stagnation assumption, and the order-1 and order-2 stability of the GWO under stagnation assumption are proved successively. In section 5, we conclude the theoretical results above.

## 2. An introduction of GWO

As mentioned above, the GWO algorithm is inspired by the leadership hierarchy and hunting mechanism of grey wolves in nature. In the GWO, the best, the second best, and the third best search agents are called α, β, and δ wolves, respectively, and other search agents are called ω wolves. Let $T$ be the total number of iterations of the GWO, $D$ be the dimension of the optimization problem, $N$ be the number of search agents in the group, and the positions of α, β, and δ wolves at iteration $t$ are marked as $\mathbf{p_1}(t) = [p_{11}(t), \dots, p_{1D}(t)] = [p_{1j}(t)|j = 1, \dots D]$, $\mathbf{p_2}(t) = [p_{2j}(t)|j = 1, \dots D]$, $\mathbf{p_3}(t) = [p_{3j}(t)|j = 1, \dots D]$, respectively. The position of the $i$th ω wolf at iteration $t$ $\mathbf{x_i}(t) = [x_{ij}(t)|i = 1, \dots, N; j = 1, \dots D]$ is updated by the guidance of the positions of α, β, and δ wolves $\mathbf{p_k}(t)(k = 1,2,3)$, which can be formulated as:

$$\mathbf{x_i}(t+1) = \frac{1}{3}\sum_{k=1}^{3} \mathbf{x'_k}(t) \tag{2.1}$$

where $\mathbf{x'_k}(t) = [x'_{kj}(t)|k = 1,2,3; j = 1, \dots D]$ is an intermediate vector, and $x'_{kj}(t)$ can be calculated as:

$$x'_{kj}(t) = p_{kj}(t) + A_{kj}|C_{kj}p_{kj}(t) - x_{ij}(t)| \tag{2.2}$$

where $A_{kj} \sim U[-a(t), a(t)]$ ($U[p, q]$ denotes the uniform distribution with lower bound $p$ and upper bound $q$, similarly hereinafter), $C_{kj} \sim U[0,2]$, and $a(t)$ is linearly decreased from 2 to 0 during iteration:

$$a(t) = 2\left(1 - \frac{t}{T}\right) \tag{2.3}$$

Obviously, $\mathbf{x'_1}(t), \mathbf{x'_2}(t),$ and $\mathbf{x'_3}(t)$ represent the guidance of the positions of α, β, and δ wolves to that of the $i$th ω wolf, respectively. According to Eq. (2.1), the new position of the $i$th search agent $\mathbf{x_i}(t+1)$ is an average of the three intermediate vectors $\mathbf{x'_k}(t)(\mathbf{k} = 1,2,3)$, which provides an intuitive reflection of the hunting process of ω wolves in real life. In addition, when the fitness of the new search agent $\mathbf{x_i}(t+1)$ is better than that of the best three search agents at iteration $t$ $\mathbf{p_k}(t)$, the new position of the best three search agents $\mathbf{p_k}(t+1)$ will be replaced by $\mathbf{x_i}(t+1)$, which can be formulated as:

$$\mathbf{p_1}(t+1) = \begin{cases} \mathbf{x_i}(t+1), f(\mathbf{x_i}(t+1)) < f(\mathbf{p_1}(t)) \\ \mathbf{p_1}(t), \text{otherwise} \end{cases} \tag{2.4}$$

$$\mathbf{p_2}(t+1) = \begin{cases} \mathbf{x_i}(t+1), f(\mathbf{p_1}(t)) \leq f(\mathbf{x_i}(t+1)) < f(\mathbf{p_2}(t)) \\ \mathbf{p_2}(t), \text{otherwise} \end{cases} \tag{2.5}$$

$$\mathbf{p_3}(t+1) = \begin{cases} \mathbf{x_i}(t+1), f(\mathbf{p_2}(t)) \leq f(\mathbf{x_i}(t+1)) < f(\mathbf{p_3}(t)) \\ \mathbf{p_3}(t), \text{otherwise} \end{cases} \tag{2.6}$$

where the optimization problem is assumed to be $\min f(\mathbf{x})$. In summary, the pseudo-code of the GWO is given in Fig. 1.

```
1:   Initialize the grey wolf population  x_i(1)(i = 1, ..., N)
2:   Calculate the fitness of each search agent
3:   Name the best search agent as  p_1(1)
4:   Name the second best search agent as  p_2(1)
5:   Name the third best search agent as  p_3(1)
6:   while  t ≤ T
7:       for each search agent
8:           Update the position of the current search agent  x_i(t + 1)  by
(2.1)-(2.3)
9:       end for
10:      Calculate the fitness of all search agents
11:      Update the position of the best three search agents  p_1(t + 1), p_2(t +
1), p_3(t + 1)  by (2.4)-(2.6)
12:      t = t + 1
13:  end while
14:  return  p_1(T)  and  f(p_1(T))
```

Fig. 1.   Pseudo-code of the GWO.

## 3. Analytical model of sampling distribution of $x_i(t+1)$ with constant $x_i(t)$ assumption

In this section, the stagnation assumption is introduced and explained under the assumption that $x_i(t)$ is treated as a constant, and the sampling distribution characteristics of the positions of the intermediate vector $x'_k(t)$ as well as the new solution $x_i(t+1)$ are derived and analyzed successively.

### 3.1 Preliminaries and assumptions

According to Eqs. (2.1) - (2.3), the position of the *i*th search agent at iteration $t+1$ (marked as $x_i(t+1)$ is a function of the position of the *i*th search agent at iteration $t$ $x_i(t)$, the position of the best three search agents at iteration $t$ $p_k(t)(k=1,2,3)$, and the random parameters $A_{kj}, C_{kj}$. Obviously, the position vectors $x_i(t)$ and $p_k(t)$ are constantly changing during iterations, thus the $x_i(t)$ should be seen as a complicated stochastic process corresponding to the iteration *t*, which brings a certain degree of difficulty to the analysis and solution of the sampling distribution characteristics of $x_i(t+1)$. Accordingly, in view of the convenience of problem analysis, it is necessary to put forward some appropriate assumptions.

In the part I of the paper, the stagnation assumption is introduced and applied to the theoretical analysis of the GWO, i.e., the positions of the best three search agents $p_k(t)$ keep unchanged during iterations (as a result, in this section and the next section, the $p_k(t)$ will be abbreviated as $p_k$). As a common assumption in the field of the analysis of meta-heuristic algorithms, the stagnation assumption has been widely applied to simplify the sampling distribution and stability analysis of the particle swarm optimizer (PSO), in which the global best (gbest) of the population and the personal best (pbest) of each search agent are set to be constant [46-47]. Obviously, in the GWO, the stagnation assumption is usually approximately satisfied at the later stage. More specifically, in the earlier stage of the GWO (please see Fig. 2(a)), since the search agents are distributed in the whole feasible region of the optimization problem, the searching space of each agent is relatively large, and the population holds a wide searching range (represented by the big cube $D_i$ as shown in Fig. 2(a)), thus the ability to search for the neighborhood of the global optimum (represented by the small cube $D_g$ and point $\mathbf{x_g}$ as shown in Fig. 2(a)) is guaranteed. However, in the later stage of the GWO, the stagnation assumption is approximately satisfied, i.e. the relationship $\mathbf{p_1} \approx \mathbf{p_2} \approx \mathbf{p_3}$ holds, and the positions of all other search agents are close to that of $\mathbf{p_k}$. In the optimistic scenarios, the $\mathbf{p_k}$ is within a small neighborhood of the global optimum, which means the GWO succeeds in finding the global optimum; however, in the pessimistic scenarios (please see Fig. 2(b)), the $\mathbf{p_k}$ is quite far from the global optimum, thus the ability of the population to search for the global optimum is no longer guaranteed.

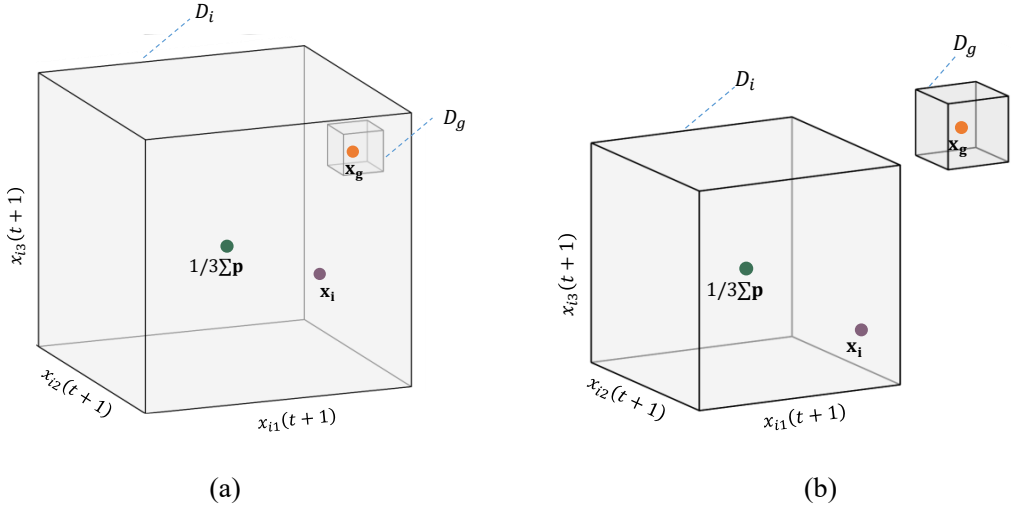

Fig. 2. Diagrams of the global optimum searching ability of the population in GWO.

With the consideration of stagnation assumption, the positions of the best three search agents can be seen independent of iteration $t$, however, the position of search agent $i$ $\mathbf{x_i}(t)$ is still a stochastic process corresponding to $t$. In this section, we further assume that the $\mathbf{x_i}(t)$ is a constant vector (thus abbreviated as $\mathbf{x_i}$ in this section), in which case the $\mathbf{x_i}(t+1)$ is no longer a stochastic process, but merely a random variable calculated as (according to Eqs. (2.1) and (2.2)):

$$x_{ij}(t+1) = \tfrac{1}{3}\sum p_j + \tfrac{1}{3}\sum_{k=1}^{3} A_{kj}|C_{kj}p_{kj} - x_{ij}| \tag{3.1}$$

where $\sum p_j = \sum_{k=1}^{3} p_{kj}$, similarly hereinafter.

According to Eq. (3.1), the $x_{ij}(t+1)$ becomes a function of six random variables $A_{kj}, C_{kj}(k = 1,2,3)$, thus the difficulty of analyzing the sampling distribution characteristics of $\mathbf{x_i}(t+1)$ is substantially reduced.

**Remark 1**

In order to better understand the work to be done in this section, a simplified solution updating equation applied to evolution strategy (ES) [48] and evolutionary programming (EP) [49] algorithm as given in (3.2) is utilized to further explain our purpose:

$$x_{ij}(t+1) = x_{ij}(t) + N(0,\sigma^2) \tag{3.2}$$

where $N(0,\sigma^2)$ denotes a random variable following a normal distribution with mean 0 and variance $\sigma^2$. Although the distribution of the stochastic process $x_{ij}(t)$ is hard to derive, we can easily conclude that the new solution $x_{ij}(t+1)$ follows a normal distribution with mean $x_{ij}(t)$ and variance $\sigma^2$ under the assumption that $x_{ij}(t)$ is treated as a constant, which points out an obvious direction for researchers to control the mutation scale of the ES/EP by adjusting the parameter $\sigma$ [50]. Furthermore, extensive research has been carried out to improve the performances of the ES and EP by changing the $N(0,\sigma^2)$ given in Eq. (3.2) into random variables following other kinds of distributions (such as Cauchy distribution [51], Levy distribution [52], etc. [53]). All of the research above show the significance of studying the sampling distribution characteristics of the new solution $x_{ij}(t+1)$ with assumption that the original solution $x_{ij}(t)$ is treated as a constant. As shown in Eqs. (2.1) and (2.2), the solution updating equation of the GWO is designed by straightly mimicking the predation behavior of grey wolves rather than providing the distribution types (or other detailed characteristics) of the updated solutions, therefore in this section we aim at analyzing and modeling the sampling distribution characteristics of the updated solution under the stagnation and constant $\mathbf{x_i}(t)$ assumptions.

According to Eq. (2.1), the $\mathbf{x_i}(t+1)$ is a linear combination of three items with the same form $\mathbf{x'_k}(t)$ (abbreviated as $\mathbf{x'_k}$ in this section), which means the joint probability density function (PDF) of $\mathbf{x_i}(t+1)$ can be derived by calculating the convolution among those of $\mathbf{x'_k}(k=1,2,3)$. Thus, we will firstly analyze and derive the sampling distribution characteristics of $\mathbf{x'_k}$ in section 3.2, and then discuss the sampling distribution characteristics of $\mathbf{x_i}(t+1)$ in section 3.3.

## 3.2 Sampling distribution analysis of $\mathbf{x'_k}$

In this section, under the stagnation and constant $\mathbf{x_i}(t)$ assumptions, the sampling distribution of $\mathbf{x'_k}$ is firstly derived and discussed. According to Eq. (3.1), each component of the new solution $x_{ij}(t+1)(j=1,...,D)$ is calculated independently. In addition, all of the $2D$ random parameters $A_{kj}, C_{kj}(j=1,...,D)$ are independent of each other, thus all of the $D$ random variables $x_{ij}(t+1)(j=1,...,D)$ are independent of each other. Consequently, let the joint PDF of $\mathbf{x'_k}$ be $c_{kj}(\mathbf{u})(\mathbf{u}=[u_j|j=1,...,D])$ and the PDF of $x'_{kj}$ be $g_{kj}(u_j)$, we have (all necessary proofs of propositions and corollaries can be found in Appendix):

**Proposition 1.1**

$$c_{kj}(\mathbf{u}) = \prod_{j=1}^{D} g_{kj}(u_j) \tag{3.3}$$

where the expression of $g_{kj}(u_j)$ can be derived as follows:

**Proposition 1.2**

Let
$$m_{kj} = a(t)(-|p_{kj}|+|x_{ij}-p_{kj}|), n_{kj} = a(t)(|p_{kj}|+|x_{ij}-p_{kj}|) \tag{3.4}$$

When $m_{kj} \leq 0$,

$$g_{kj}(u_j) = \begin{cases} \frac{1}{n_{kj}-m_{kj}} \ln \frac{\sqrt{-m_{kj}n_{kj}}}{|u_j-p_{kj}|}, & |u_j-p_{kj}| < -m_{kj} \\ \frac{1}{2(n_{kj}-m_{kj})} \ln \frac{n_{kj}}{|u_j-p_{kj}|}, & -m_{kj} \leq |u_j-p_{kj}| < n_{kj} \\ 0, & |u_j-p_{kj}| \geq n_{kj} \end{cases} \tag{3.5}$$

Otherwise,

$$g_{kj}(u_j) = \begin{cases} \frac{1}{2(n_{kj}-m_{kj})} \ln \frac{n_{kj}}{m_{kj}}, & |u_j - p_{kj}| < m_{kj} \\ \frac{1}{2(n_{kj}-m_{kj})} \ln \frac{n_{kj}}{|u_j-p_{kj}|}, & m_{kj} \leq |u_j - p_{kj}| < n_{kj} \\ 0, & |u_j - p_{kj}| \geq n_{kj} \end{cases} \quad (3.6)$$

As the core conclusion of this section, the proposition 1.2 gives an exact analytical expression for the PDF of intermediate variable $x'_{kj}$. According to the proposition 1.2, obviously we have:

**Corollary 1.1**

The domain of $g_{kj}(u_j)$ is $\{u_j \in \mathbb{R} | g_{kj}(u_j) > 0\} = (p_{kj} - n_{kj}, p_{kj} + n_{kj})$ (in this corollary, the domain is defined as $\{u_j \in \mathbb{R} | g_{kj}(u_j) > 0\}$, however, in the rest of the paper, for the sake of simplicity, we might say that a certain function $h$ holds some characteristics within $\mathbb{R}$, yet the domain of $h$ as defined in this corollary is only a subset of $\mathbb{R}$, and the value of $h$ that exceeds the domain (as defined in this corollary) is equal to 0). The curve of $g_{kj}(u_j)$ is symmetrical about $u_j = p_{kj}$, and $g_{kj}(u_j)$ is non-decreasing within $(p_{kj} - n_{kj}, p_{kj})$.

From corollary 1.1, the distribution of $x'_{kj}$ is a bounded single-peak distribution symmetrical about $p_{kj}$.

According to the proposition 1.2, the curves of $g_{kj}(u_j)$ under cases that $m_{kj} \leq 0$ and $m_{kj} > 0$ can be drawn as Fig. 3 (a) and Fig. 3 (b), respectively, in which all curve parts correspond to the logarithmic functions. It is worth noting that when $m_{kj} \leq 0$, the curve of $g_{kj}(u_j)$ within $(p_{kj} - n_{kj}, p_{kj})$ is formed by two different logarithmic curves (i.e. $g_{kj}(u_j)$ is not differentiable at $u_j = p_{kj} \pm m_{kj}$). Furthermore, when $m_{kj} \leq 0$, the range of $g_{kj}(u_j)$ is not bounded, and $g_{kj}(p_{kj}) = +\infty$ holds.

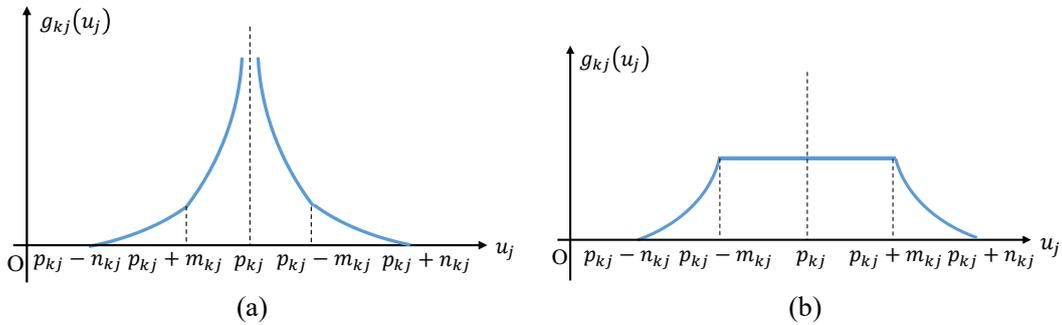

Fig. 3. Curves of $g_{kj}(u_j)$, (a): $m_{kj} \leq 0$, (b): $m_{kj} > 0$.

The proposition 1.2 can be verified by Monte Carlo (MC) simulations. Fig. 4 shows the frequency histogram of $x'_{kj}$ obtained by MC simulation and its PDF curve given in Eqs. (3.4) – (3.6) with four sets of values of variables $a(t), p_{kj}, x_{ij}$, and $m_{kj}$. For each set of variable values, the sample size is set to be $10^6$, and the number of frequency histogram segments is set to be 60. Obviously, as the sample size and the number of frequency histogram segments increase, the frequency histogram of $x'_{kj}$ obtained by MC simulation tends to coincide with the curve of the PDF of $x'_{kj}$ (i.e. $g_{kj}(u_j)$). It can be seen from Fig. 4 that whatever the sign of $m_{kj}$ is, the frequency histogram of $x'_{kj}$ coincides with the analytical curve of $g_{kj}(u_j)$ (as calculated by Eqs. (3.4)-(3.6)) very well.

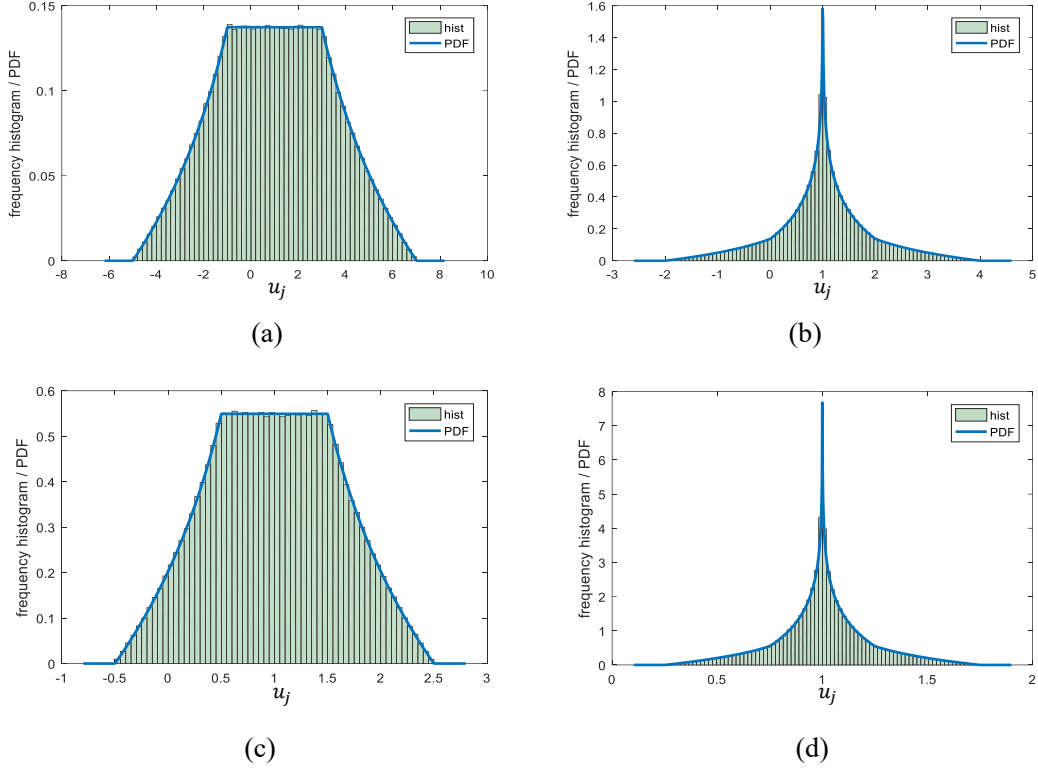

Fig. 4. Frequency histograms of $x'_{kj}$ and curves of $g_{kj}(u_j)$. (a): $a(t) = 2, p_{kj} = 1, x_{ij} = 3, m_{kj} = 2 > 0$, (b): $a(t) = 2, p_{kj} = 1, x_{ij} = 0.5, m_{kj} = -1 < 0$, (c): $a(t) = 0.5, p_{kj} = 1, x_{ij} = 3, m_{kj} = 0.5 > 0$, (d): $a(t) = 0.5, p_{kj} = 1, x_{ij} = 0.5, m_{kj} = -0.25 < 0$

Based on the analytical model of $g_{kj}(u_j)$ discussed above, in the last part of this section, the distribution characteristics of $\mathbf{x'_k}$ will be discussed elaborately. Since all of the $D$ random variables $x_{ij}(t+1)(j = 1, ..., D)$ are independent of each other, and according to the corollary 1.1, the domain of $g_{kj}(u_j)$ is $(p_{kj} - n_{kj}, p_{kj} + n_{kj})$, obviously we have:

**Corollary 1.2**

The random vector $\mathbf{x'_k}$ is distributed in the $D$-dim box:
$$D_k = \prod_{j=1}^{D}(p_{kj} - n_{kj}, p_{kj} + n_{kj}) \tag{3.7}$$

According to the corollary 1.2, the center of $D_k$ is $\mathbf{p_k}$, and the length of the $j$th side of $D_k$ is $2n_{kj}$. Obviously, the measurement of $D_k$ is an increasing function of $a(t)$ (please see Eq. (3.4)), which means the searching space of agent $i$ decreases during iterations.

The MC simulation method is applied to illustrate the distribution characteristics of $\mathbf{x'_k}$. Here, the sample size of $\mathbf{x'_k}$ is set to be 5000, and $a = 2, \mathbf{p_k} = [1,1], \mathbf{x_i} = [1.5, 3]$, then the MC simulation results of $\mathbf{x'_k}$ (marked by blue points) and the theoretical value of $D_k$ (as calculated via Eq. (3.7), marked by blue dashed rectangle) are given in Fig. 5. It can be seen that $D_k$ is approximately an envelope of $\mathbf{x'_k}$, whose center is $\mathbf{p_k}$, and the side of length is directly proportional to $a(t)$ (in Fig. 5 (a), $a(t) = 2$, the size of $D_k$ is $6 \times 12$; in Fig. 5 (b), $a(t) = 0.5$, the size of $D_k$ is $1.5 \times 3$). Furthermore, the probability that $\mathbf{x'_k}$ appears next to $\mathbf{p_k}$ is obviously larger than that next to the boundary of $D_k$, which coincides with the conclusion that $g_{kj}(u_j)$ is a single-peak function (please see the corollary 1.1).

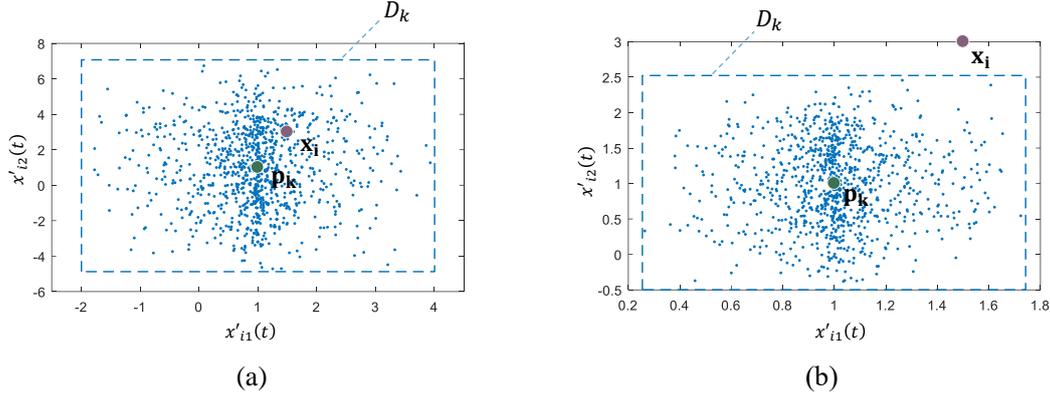

Fig. 5. MC simulation results of $\mathbf{x'_k}$ (marked by blue points), (a): $a(t) = 2$, (b): $a(t) = 0.5$

### 3.3 Analytical modeling of sampling distribution of $\mathbf{x_i}(t+1)$

In section 3.2, we derived the analytical expression and discussed the distribution characteristics that $\mathbf{x'_k}$ follows. Based on the results obtained above, in this section, we will study the distribution model that $\mathbf{x_i}(t+1)$ follows under stagnation and constant $\mathbf{x_i}(t)$ assumptions. According to Eq. (2.1), the $x_{ij}(t+1)$ is a linear combination of $x'_{kj}(k = 1,2,3)$, thus the PDF of $x_{ij}(t+1)$ can be derived by calculating the convolution among $g_{kj}(u_j)(k = 1,2,3)$ [54]. Let the PDF of $x_{ij}(t+1)$ be $h_{ij}(u_j)$, obviously we have:

**Proposition 2.1**

$$h_{ij}(u_j) = (g_{1j} * g_{2j} * g_{3j})(3u_j) \qquad (3.8)$$

Although the expression of $h_{ij}(u_j)$ is given in proposition 2.1 in the convolution form, it is almost impossible to substitute Eqs. (3.5) and (3.6) into Eq. (3.8) to derive the exact analytical expression of $h_{ij}(u_j)$. This is because the expression of $g_{kj}(u_j)$ is a piecewise function, hence the exact analytical expression of $h_{ij}(u_j)$ becomes a more complicated piecewise function with more subsections after convolution, and the endpoints of each segment require lengthy and complicated discussions, let alone the difficulty of calculating the exact analytical expression on each segment.

**Remark 2**

In order to better illustrate the complexity of the exact analytical expression of $h_{ij}(u_j)$, we try to derive the expression of $(g_{1j} * g_{2j})(u_j)$ on its leftmost subsection region. We suppose that $m_{1j}, m_{2j} > 0$, and let $g_{10}(u_j) = g_{1j}(u_j + p_{1j}), g_{20}(u_j) = g_{2j}(u_j + p_{2j})$, then both $g_{10}(u_j), g_{20}(u_j)$ are piecewise functions, whose subsection region endpoints are (from left to right) $-n_{1j}, -m_{1j}, m_{1j}, n_{1j}$ and $-n_{2j}, -m_{2j}, m_{2j}, n_{2j}$, respectively, which means the subsection region endpoints of $g_{10} * g_{20}$ are $\pm n_{1j} \pm n_{2j}, \pm m_{1j} \pm n_{2j}, \pm n_{1j} \pm m_{2j}, \pm m_{1j} \pm m_{2j}$ (there may be duplicates). The number of subsection regions may be as many as 15, let alone the tedious discussion of the relationship among these endpoints.

In addition, we try to derive the expression of the convolution results on its leftmost subsection region. Suppose that $-n_{1j} - m_{2j} < -m_{1j} - n_{2j}$, in which case the leftmost subsection region of $(g_{1j} * g_{2j})(u_j)$ becomes $(-n_{1j} - n_{2j}, -n_{1j} - m_{2j})$, and we have:

$$(g_{1j} * g_{2j})(u_j + p_{1j} + p_{2j})$$

$$= (g_{10} * g_{20})(u_j)$$

$$= \int_{-n_{1j}}^{u_j+n_{2j}} \frac{1}{2(n_{1j}-m_{1j})} \cdot \ln\frac{n_{1j}}{p_{1j}-v} \cdot \frac{1}{2(n_{2j}-m_{2j})} \cdot \ln\frac{n_{2j}}{p_{2j}-v} dv$$

$$= \frac{1}{4(n_{1j}-m_{1j})(n_{2j}-m_{2j})} \int_{-n_{1j}}^{u_j+n_{2j}} \ln\frac{n_{1j}}{p_{1j}-v} \cdot \ln\frac{n_{2i}}{p_{2j}-v} dv$$

$$= \frac{1}{4(n_{1j}-m_{1j})(n_{2j}-m_{2j})} \Big\{ v - (v-p_{1j})[-1+\ln(v-p_{1j})] + (v-p_{1j})\Big(1+$$

$$\ln\frac{n_{1j}}{v-p_{1j}}\Big)\ln\frac{n_{2j}}{v-p_{2j}} + (p_{2j}-p_{1j})\ln(v-p_{2j}) + (\ln p_{1j})[v-p_{1j}+(p_{2j}-p_{1j})\ln(v-$$

$$p_{2j})] - (p_{2j}-p_{1j})\Big[\ln(v-p_{1j})\ln\frac{v-p_{2j}}{p_{1j}-p_{2j}} + Li_2\Big(-\frac{v-p_{2j}}{p_{1j}-p_{2j}}\Big)\Big]\Big\}\Big|_{-(u_j+n_{2j})}^{n_{1j}}$$

where $Li_n(z) = \sum_{k=1}^{\infty} \frac{z^k}{k^n}$ denotes the poly-logarithm function [55]. Evidently, such a lengthy result has little practical meaning and application value. Taking the tedious classified discussion of subsection region endpoints into account, we could easily conclude that the exact analytical expression of $h_{ij}(u_j)$ becomes too lengthy to hold any practical value.

Fortunately, based on the results provided in section 3.2, we can further analyze the characteristics of $h_{ij}(u_j)$ as follows. If setting

$$\sum p_j = p_{1j} + p_{2j} + p_{3j}, \sum m_j = m_{1j} + m_{2j} + m_{3j}, \sum n_j = n_{1j} + n_{2j} + n_{3j} \quad (3.9)$$

for the sake of simplicity, then we have:

**Corollary 2.1**

$h_{ij}(u_j)$ is a continuous bounded function defined within $\{u_j \in R | h_{ij}(u_j) > 0\} = \left(\frac{1}{3}\sum p_j - \frac{1}{3}\sum n_j, \frac{1}{3}\sum p_j + \frac{1}{3}\sum n_j\right)$.

According to the corollary 2.1, the domain and range of $h_{ij}(u_j)$ can be obtained. Similar with that of $g_{kj}(u_j)$, the domain of $h_{ij}(u_j)$ is an open interval with finite length, which means when $x_{ij}(t)$ is regarded as constant, the searching space of $x_{ij}(t+1)$ is limited within a range symmetrical about $\frac{1}{3}\sum p_j$ (different from the mutation operator of the ES/EP as shown in Remark 1, if the mutation operator of ES/EP is applied to update $x_{ij}$, the searching space of $x_{ij}(t+1)$ will become the whole feasible region). In addition, it is worth noting that the non-continuous function $g_{kj}(u_j)$ with unbounded range (when $m_{kj} \leq 0$) is transformed into a continuous function $h_{ij}(u_j)$ with bounded range after the application of convolution operator, in which process the smoothing effect of convolution operator [56-57] is demonstrated, which will be elaborated in corollary 2.3.

**Corollary 2.2**

The curve of $h_{ij}(u_j)$ is symmetrical about $u_j = \frac{1}{3}\sum p_j$, $h_{ij}\left(\frac{1}{3}\sum p_j - \frac{1}{3}\sum n_j\right) = 0$, and $h_{ij}(u_j)$ is non-decreasing within $\left(\frac{1}{3}\sum p_j - \frac{1}{3}\sum n_j, \frac{1}{3}\sum p_j\right)$.

According to the corollary 2.2, after convolution, the symmetry and unimodal characteristics of $g_{kj}(u_j)$ are both retained, i.e. the $h_{ij}(u_j)$ is a single-peak function symmetrical about $u_j = \frac{1}{3}\sum p_j = \frac{1}{3}\sum_{k=1}^{3} p_{kj}$, which is exactly the average of those of $g_{kj}(u_j)(k=1,2,3)$ (i.e. $p_{kj}$).

**Corollary 2.3**

When there exists at least one positive number among $m_{kj}(k=1,2,3)$, the $h_{ij}(u_j)$ is derivable within $\mathbb{R}$, and $g'_{ij}\left(\frac{1}{3}\sum p_j\right) = 0$ holds.

The corollary 2.3 provides a further indication of the smoothing effect of convolution operator [56-57]. According to the proposition 1.2 (as shown in Fig. 3), the $g_{kj}(u_j)$ is non-differentiable at $u_j = p_{kj} \pm m_{kj}$, and when $m_{kj} \leq 0$, $g_{kj}(p_{kj}) = +\infty$, the $g_{kj}(u_j)$ is non-continuous at $u_j = p_{kj}$. However, when there exists at least one positive number among $m_{kj}(k=1,2,3)$, the smoothing effect of convolution operator makes $h_{ij}(u_j)$ a differentiable function within $\mathbb{R}$. However, the smoothing effect is limited within conditions that $\exists k, m_{kj} > 0$. When $\forall k = 1,2,3, m_{kj} \leq 0$, $g_{kj}(u_j)(k=1,2,3)$ are all non-continuous functions, then the differentiability of $h_{ij}(u_j)$ at $u_j = \frac{1}{3}\sum p_j$ is no more guaranteed.

The conclusions above can be verified by MC simulations. Here, the sample size is set to be $10^6$, and the number of frequency histogram segments is set to be 60. Since the sign of $m_{kj}$ has a significant impact on the shape of the curves of $g_{kj}(u_j)$ and $h_{ij}(u_j)$, the simulations are carried out under the following four cases as shown in Eq. (3.10). In addition, in each case two groups of values of $x_{ij}$ and $p_{kj}$ are randomly chosen, thus 8 groups of the frequency histograms of $x_{ij}(t+1)$ are generated and shown in Fig. 6. In order to better exhibit the characteristics of $h_{ij}(u_j)$ generated by MC simulations, the upper-side midpoints of all rectangles in each frequency histogram are connected together as the curve of $h_{ij}(u_j)$. For reader's reference, all parameters corresponding to each frequency histogram are given in Table I.

$$\sum_{k=1}^{3} \varepsilon(m_{kj}) = 0,1,2,3, \text{ where } \varepsilon(x) = \begin{cases} 1, x > 0 \\ 0, x \leq 0 \end{cases} \quad (3.10)$$

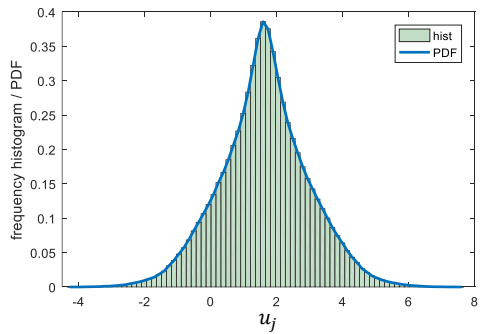

(a)

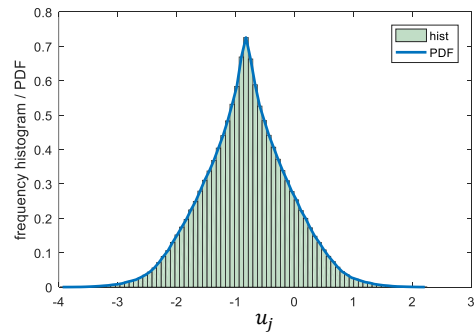

(b)

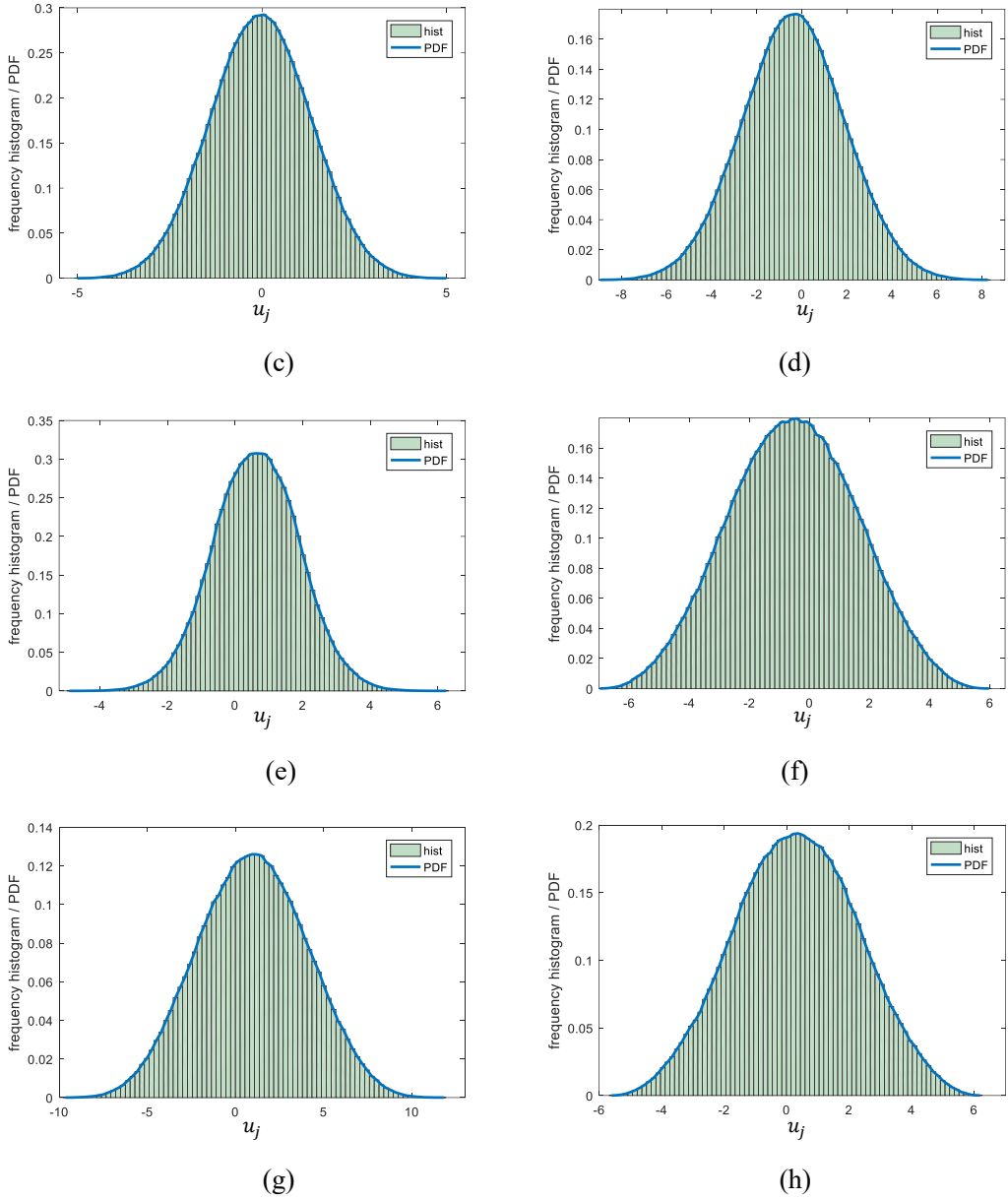

(c) (d) (e) (f) (g) (h)

Fig. 6. Frequency histograms, please see Table I for parameters corresponding to each sub-graph.

Table I Parameters corresponding to each sub-graph as shown in Fig. 6.

| Subfigure No. | $p_{1j}$ | $p_{2j}$ | $p_{3j}$ | $\frac{1}{3}\sum p_j$ | $x_{ij}$ | $\sum_{k=1}^{3} \varepsilon(m_{kj})$ |
|---|---|---|---|---|---|---|
| a | 0.43 | 1.74 | 2.70 | 1.62 | 0.83 | 0 |
| b | -0.18 | -0.97 | -1.31 | -0.82 | -0.22 | 0 |
| c | 0.57 | 0.86 | -1.56 | -0.04 | 0.72 | 1 |
| d | 0.65 | 1.13 | -2.89 | -0.37 | 0.90 | 1 |
| e | 0.09 | -1.07 | 2.88 | 0.63 | 2.08 | 2 |
| f | 0.56 | -1.85 | -0.34 | -0.54 | -3.14 | 2 |
| g | -0.30 | 0.96 | 2.24 | 0.97 | -2.87 | 3 |
| h | 0.22 | 0.55 | 0.18 | 0.32 | 3.11 | 3 |

It can be seen from Fig. 6 that whatever the values of $x_{ij}$ and $p_{kj}$ are, the $h_{ij}(u_j)$ are always single-peak functions symmetrical about $u_j = \frac{1}{3}\sum p_j$. In addition, when there exists at least one positive number among $m_{kj}(k = 1,2,3)$ (corresponding to Fig. 6 (c) – (h)), the curve of $h_{ij}(u_j)$ is smooth, and the $h_{ij}(u_j)$ is derivable within $\mathbb{R}$; when $\forall k = 1,2,3, m_{kj} \leq 0$ holds (corresponding to Fig. 6 (a) – (b)), the curve of $h_{ij}(u_j)$ is not guaranteed to be smooth at $u = \frac{1}{3}\sum p_j$, which means the $h_{ij}(u_j)$ might not be derivable at $u_j = \frac{1}{3}\sum p_j$. Thus, corollaries 2.1 – 2.3 are verified.

In the last part of this section, the distribution characteristics of $\mathbf{x_i}(t + 1)$ is discussed. Since all of the $D$ random variables $x_{ij}(t+1)(j = 1, \ldots, D)$ are independent of each other, let the joint PDF of $\mathbf{x_i}(t + 1)$ be $c_i(\mathbf{u})$, then we have:

**Proposition 2.2**

$$c_i(\mathbf{u}) = \prod_{j=1}^{D} h_{ij}(u_j) \tag{3.11}$$

Based on the analytical expression and characteristics of $g_{ij}(u_j)$ discussed above, we have:

**Corollary 2.4**

The random vector $\mathbf{x_i}(t + 1)$ is distributed in the $D$-dim box:

$$D_i = \prod_{j=1}^{D}\left(\frac{1}{3}\sum p_j - \frac{1}{3}\sum n_j, \frac{1}{3}\sum p_j + \frac{1}{3}\sum n_j\right) \tag{3.12}$$

According to the corollary 2.4, the center of $D_i$ is $\frac{1}{3}\sum \mathbf{p}$, and the length of the $j$th side of $D_i$ is $\frac{2}{3}\sum n_j$. Similar to the characteristics of $D_k$ (as given in the corollary 1.2), the measurement of $D_i$ is an increasing function of $a(t)$ (please see Eq. (3.9)), which means the searching space of agent $i$ decreases during iterations. Here, the MC simulation is applied to illustrate the distribution of $\mathbf{x_i}(t + 1)$. The sample size of $\mathbf{x_i}(t + 1)$ is set to be 5000, and $D = 2, \mathbf{p_1} = [0,0], \mathbf{p_2} = [0,1], \mathbf{p_3} = [1,0], \mathbf{x_i} = [1.5,3]$, then the MC simulation results of $\mathbf{x_i}(t + 1)$ (marked by blue points) and the theoretical values of $D_i$ (as calculated via Eq. (3.12), and marked by blue dotted rectangles) are shown in Fig. 7.

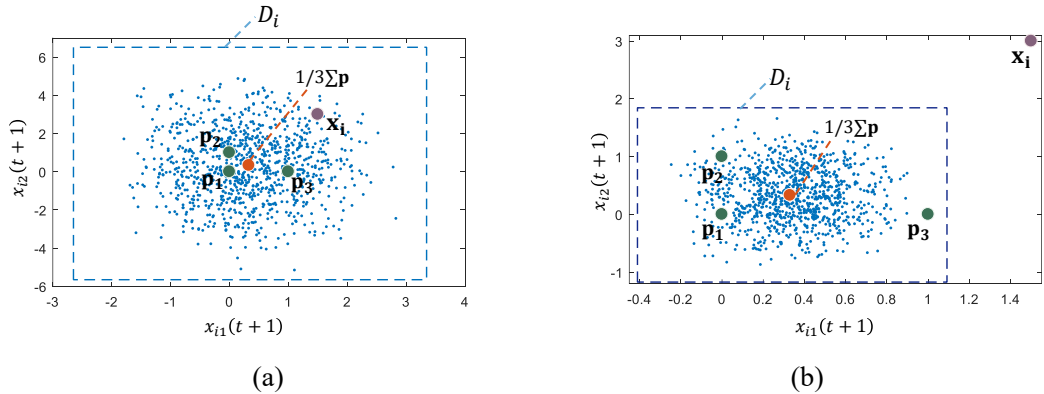

Fig. 7. MC simulation results of $\mathbf{x_i}(t + 1)$ (blue points), (a): $a(t) = 2$, (b): $a(t) = 0.5$

It can be seen from Fig. 7 that the distribution of $\mathbf{x_i}(t + 1)$ approximately forms a rectangle, whose center is $\frac{1}{3}\sum \mathbf{p}$ (marked by the red dot in Fig. 7). In addition, the probability that $\mathbf{x_i}(t + 1)$ appears next to $\frac{1}{3}\sum \mathbf{p}$ is obviously larger than the probability of being far away from $\frac{1}{3}\sum \mathbf{p}$, which coincides well with the conclusion that the $h_{ij}(u_j)$ is a single-peak function (please see the corollary 2.2). However, different from the MC simulation results of $\mathbf{x'_k}$, the side-length of the rectangle envelope of the simulation results of $\mathbf{x_i}(t+1)$ (see the blue points as shown in Fig. 7) is evidently smaller than that of the theoretical value of $D_i$ (calculated by Eq. (3.12), and marked by the blue dotted rectangle as shown in Fig. 7), which indicates that the probability that $x_{ij}(t+1)$ becomes close to its theoretical maximum (or minimum) value is remarkably smaller than that of $x_{kj}(t)$.

In order to illustrate this phenomenon more directly, we generate the curves of $h_{ij}(u_j)$ and $g_{kj}(u_j)$ in the same coordinate system by MC simulations as shown in Fig. 8 (a). Let $a(t) = 2, p_{1j} = p_{2j} = p_{3j} = 1, x_{ij} = 3$, then it can be calculated that both the domains of $h_{ij}(u_j)$ and $g_{kj}(u_j)$ are $[-5,7]$ (please see the corollary 1.1 and corollary 2.1). In addition, based on the numerical results of $h_{ij}(u_j)$ and $g_{kj}(u_j)$, the CDF of $x'_{kj}(t)$ (marked as $G_{kj}(u_j)$) and the CDF of $x_{ij}(t+1)$ (marked as $H_{ij}(u_j)$) are calculated via Trapezoid formula [58], and the curves of $G_{kj}(u_j)$ and $H_{ij}(u_j)$ are drawn in the same coordinate system as well, which can be found in Fig. 8 (b). It can be seen that when $u_j$ is close to -5, the value of $H_{ij}(u_j)$ is obviously smaller than that of $G_{kj}(u_j)$, and the value of $h_{ij}(u_j)$ is obviously smaller than that of $g_{kj}(u_j)$, which indicates that the probability that $\mathbf{x_i}(t+1)$ becomes close to its boundary is remarkably smaller than that of $\mathbf{x'_k}$.

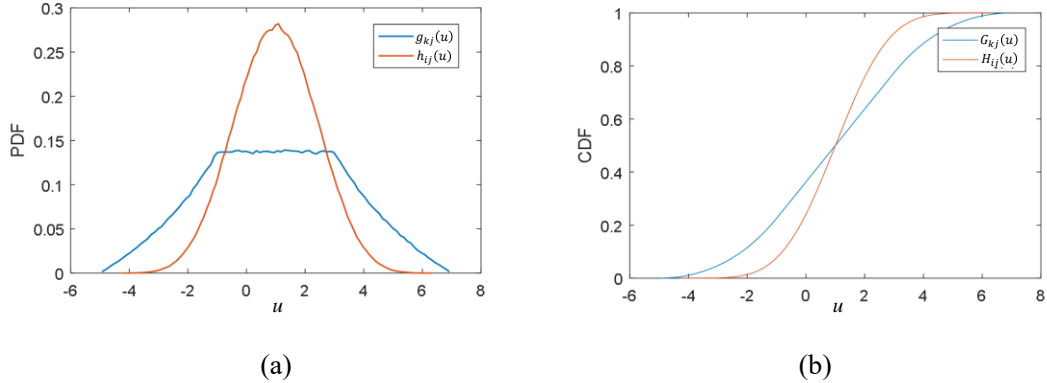

(a)          (b)

Fig. 8. (a): Curves of $h_{ij}(u_j)$ and $g_{kj}(u_j)$, (b): Curves of $H_{ij}(u_j)$ and $G_{kj}(u_j)$.

## 4. Analytical model of sampling distribution of $\mathbf{x_i}(t+1)$ without constant $\mathbf{x_i}(t)$ assumption

In section 3, under both stagnation and constant $\mathbf{x_i}(t)$ assumptions, we derived and analyzed the distribution that $\mathbf{x_i}(t+1)$ follows. Based on these works, in this section, only the stagnation assumption is considered to analyze the sampling distribution that $\mathbf{x_i}(t+1)$ follows, in which condition $\mathbf{x_i}(t+1)$ should be regarded as a stochastic process rather than a simple random variable. Although the exact analytical expression of the joint PDF of $\mathbf{x_i}(t)$ is hard to derive, several basic characteristics of the distribution can be obtained. Firstly, we will discuss the shape of the PDF of $x_{ij}(t)$, and then derive the moments and central moments of $x_{ij}(t)$, which will serve as important tools to understand the distribution that $x_{ij}(t)$ follows. Finally, as an application of the results derived above, the order-1 and order-2 stability of GWO under the stagnation assumption is introduced and proved.

## 4.1 Characteristics of the PDF of $\mathbf{x_i}(t)$

In this section, in order to better illustrate that the distribution of $x_{ij}(t)$ is related to iteration $t$, we mark the PDF of $x_{ij}(t+1)$ as $h_{ij}(u_j, t+1)$ rather than $h_{ij}(u_j)$ (as marked in section 3).

The domain of $h_{ij}(u_j, t)$ is firstly discussed. Owing to the existence of absolute value sign of the domain of $h_{ij}(u_j)$ under constant $x_{ij}(t)$ assumption (please see corollary 2.1, Eqs. (3.4) and (3.9)), it is hard to derive the analytical expression of the domain of $h_{ij}(u_j, t)$ recursively. However, we could easily prove that:

**Proposition 3.1**

When $t \geq 2$, the domain of $h_{ij}(u_j, t)$ contains $\left(\frac{1}{3}\sum p_j - d_1 \prod_{\tau=1}^{t-1} a(\tau), \frac{1}{3}\sum p_j + d_1 \prod_{\tau=1}^{t-1} a(\tau)\right)$, where $d_1 = \max_{x_{ij}(1)} \left|x_{ij}(1) - \frac{1}{3}\sum p_j\right|$.

Proposition 3.1 shows that the domain of $h_{ij}(u_j, t)$ contains an open interval with center $\frac{1}{3}\sum p_j$ and length $2d_1 \prod_{\tau=1}^{t-1} a(\tau)$. It can be seen that when $t < \frac{T}{2}, a(t) > 1$, the length is an increasing function of iteration $t$, which means in the first half period of the GWO, the searching range of each agent increases.

The next conclusion is pertinent to the shape of the curve of $h_{ij}(u_j, t)$:

**Proposition 3.2**

$\forall t \geq 2, h_{ij}(u_j, t)$ is a single-peak function symmetrical about $u_j = \frac{1}{3}\sum p_j$.

According to the proposition 3.2, even if the assumption that $\mathbf{x_i}(t)$ is constant is eliminated, the PDF of $x_{ij}(t+1)$ is still a single-peak function symmetrical about $u_j = \frac{1}{3}\sum p_j$ (please see the corollary 2.2 for the conclusion under constant $\mathbf{x_i}(t)$ assumption).

The propositions above can be verified by MC simulations. Here, the sample size is set to be $10^6$, and the number of frequency histogram segments is set to be 80. Let $p_{1j} = -1, p_{2j} = 1.5, p_{3j} = 2.5$ (i.e. $\frac{1}{3}\sum p_j = 1$), $x_{ij}(1) \sim U[-4,4]$ (corresponding to the optimization problem with feasible solution $[-4,4]$), maximum iterations $T = 60$ (i.e. $a(t) = 2\left(1 - \frac{t}{60}\right)$), and iteration $t = 2, 5, 10, 20, 30, 50$. The corresponding obtained frequency histograms of $x_{ij}(t)$ are shown in Fig. 9 (a) – (f), respectively.

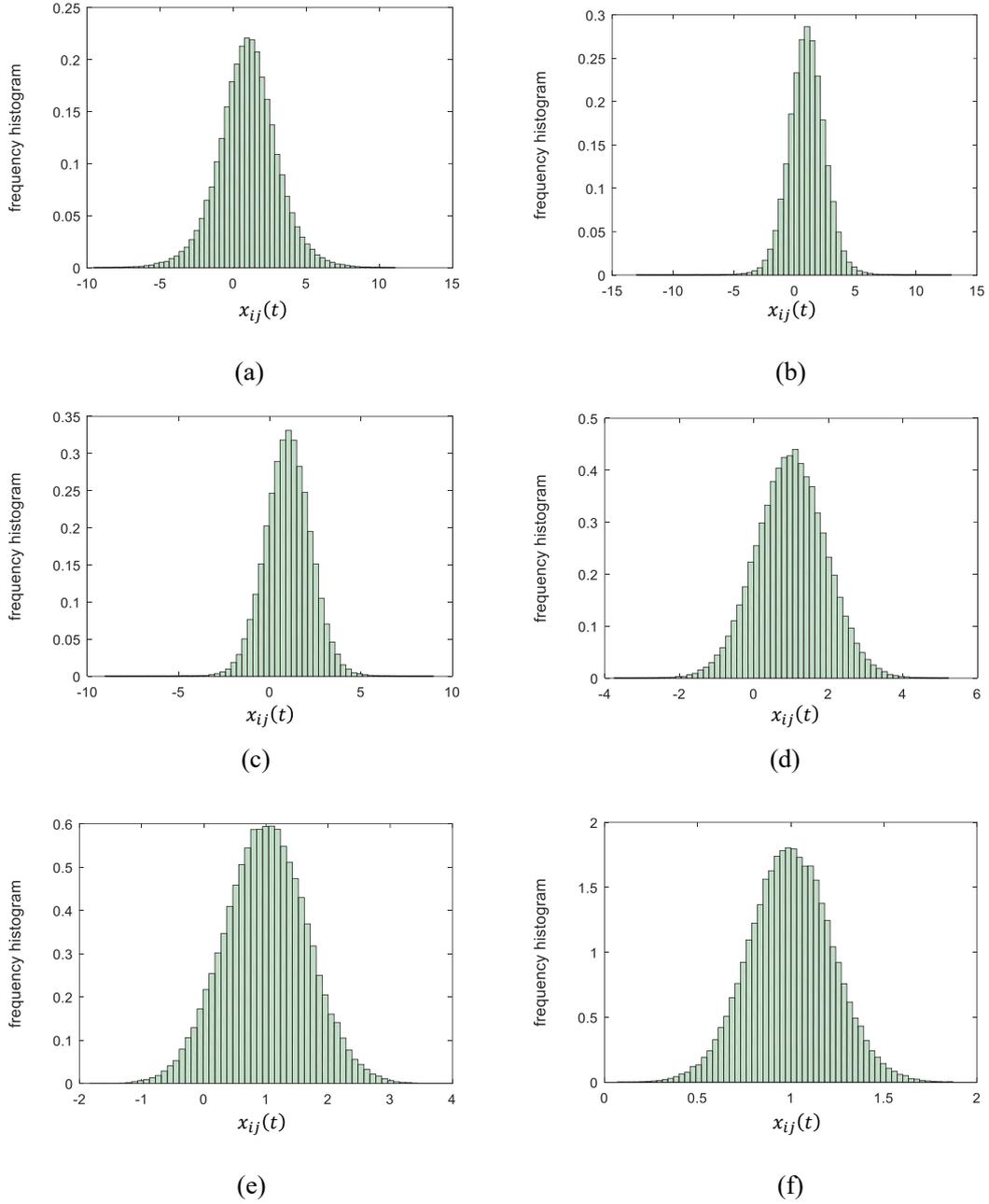

Fig. 9. Frequency histogram of $x_{ij}(t)$ under stagnation assumption, (a)-(f): $t = 2, 5, 10, 20, 30, 50$.

It can be seen from Fig. 9 that whatever the value of $t$ is, the curve of the PDF of $x_{ij}(t)$ is unimodal and symmetrical about $u_j = \frac{1}{3}\sum p_j = 1$. However, as iteration *t* increases, the distribution of $x_{ij}(t)$ continuously becomes more close to $\frac{1}{3}\sum p_j$ (please see the variation of scales of X-axis from Fig. 9 (a) to Fig. 9 (f)), which seems to be contradictory with proposition 3.1. This is because the probability that $x_{ij}(t)$ is close to its bound of domain is rather small (please see Fig. 8 and related discussions), and the probability becomes even smaller after iterations. Thus, although the domain of $h_{ij}(u_j, t)$ expands during the first half period of the GWO, it is safe to assume that the distribution of $x_{ij}(t)$ becomes closer to its center $\frac{1}{3}\sum p_j$ during iterations.

The proposition 3.1 provides a preliminary study about the characteristics of $h_{ij}(u_j,t)$. In order to better understand the distribution that $x_{ij}(t)$ follows, in the next section, we will derive the recursive expression of the central moments of any positive integer order of $x_{ij}(t)$.

**4.2 Moments of $x_{ij}(t)$**

It is known that the moment is a crucial characteristic of a random variable. In fact, if the moment of any order of a random variable is known, then the distribution of the random variable will be defined uniquely. As for the $x_{ij}(t)$, although the PDF of $x_{ij}(t)$ cannot be derived directly, the moments of any positive integer order of $x_{ij}(t)$ can be calculated as follows.

Firstly, according to the proposition 3.2, $\forall t, \mathbb{E} x_{ij}(t) = \frac{1}{3}\sum p_j$ holds. Here, we mark the $r$-order central moment of random variable $X$ as $\sigma^r X$, according to Eqs. (2.1) and (2.2), under the stagnation assumption, the $r$-order central moment of $x_{ij}(t+1)$ can be calculated as $\sigma^r x_{ij}(t+1) = \mathbb{E}\left(x_{ij}(t+1) - \mathbb{E}x_{ij}(t+1)\right)^r = \mathbb{E}\left(x_{ij}(t+1) - \frac{1}{3}\sum p_j\right)^r$, i.e.

$$\sigma^r x_{ij}(t+1) = \frac{1}{3^r}\mathbb{E}\left(\sum_{k=1}^{3} A_{kj}|C_{kj}p_{kj} - x_{ij}(t)|\right)^r \tag{4.1}$$

Since $A_{kj} \sim U[-a(t), a(t)]$, $C_{kj} \sim U[0,2]$, for a positive integer $r$, it can be calculated that $\mathbb{E}A_{kj}^r = \begin{cases} \frac{a(t)^r}{r+1}, & r \text{ is even} \\ 0, & r \text{ is odd} \end{cases}$ and $\mathbb{E}C_{kj}^r = \frac{2^r}{r+1}$. Under the stagnation assumption, $x_{ij}(t), A_{kj}, C_{kj}, p_{kj}$ are obviously independent of each other, thus Eq. (4.1) can be expanded as:

$\sigma^r x_{ij}(t+1)$

$= \frac{1}{3^r}\sum_{p=0}^{\frac{r}{2}}\sum_{q=0}^{\frac{r}{2}-p}\binom{r}{2p}\binom{r-2p}{2q}\mathbb{E}\left(A_{1j}^{2p}|C_{1j}p_{1j} - x_{ij}(t)|^{2p} A_{2j}^{2q}|C_{2j}p_{2j} - x_{ij}(t)|^{2q} A_{3j}^{2s}|C_{3j}p_{3j} - x_{ij}(t)|^{2s}\right)$

$= \frac{a(t)^r}{3^r}\sum_{p=0}^{\frac{r}{2}}\sum_{q=0}^{\frac{r}{2}-p}\frac{\binom{r}{2p}\binom{r-2p}{2q}}{(2p+1)(2q+1)(2s+1)}\mathbb{E}\left(\sum_{l=0}^{2p}\binom{2p}{l}C_{1j}^{2p-l}p_{1j}^{2p-l}x_{ij}^{l}(t)\right)$

$\left(\sum_{l=0}^{2q}\binom{2q}{l}C_{2j}^{2q-l}p_{2j}^{2q-l}x_{ij}^{l}(t)\right)\left(\sum_{l=0}^{2s}\binom{2s}{l}C_{3j}^{2s-l}p_{3j}^{2s-l}x_{ij}^{l}(t)\right)$

$= \frac{a(t)^r}{3^r}\sum_{p=0}^{\frac{r}{2}}\sum_{q=0}^{\frac{r}{2}-p}\frac{\binom{r}{2p}\binom{r-2p}{2q}}{(2p+1)(2q+1)(2s+1)}\mathbb{E}(\sum_{l_1=0}^{2p}\sum_{l_2=0}^{2q}\sum_{l_3=0}^{2s}\binom{2p}{l}\binom{2q}{l}\binom{2s}{l}p_{1j}^{2p-l_1}p_{2j}^{2q-l_2}p_{3j}^{2s-l_3}$

$C_{1j}^{2p-l_1}C_{2j}^{2q-l_2}C_{3j}^{2s-l_3}x_{ij}^{l_1+l_2+l_3}(t))$

$= \frac{a(t)^r}{3^r}\sum_{p=0}^{\frac{r}{2}}\sum_{q=0}^{\frac{r}{2}-p}\frac{\binom{r}{2p}\binom{r-2p}{2q}}{(2p+1)(2q+1)(2s+1)}(\sum_{l_1=0}^{2p}\sum_{l_2=0}^{2q}\sum_{l_3=0}^{2s}\frac{2^{r-(l_1+l_2+l_3)}\binom{2p}{l}\binom{2q}{l}\binom{2s}{l}}{(2p-l_1+1)(2q-l_2+1)(2s-l_3+1)}$

$p_{1j}^{2p-l_1}p_{2j}^{2q-l_2}p_{3j}^{2s-l_3}\mathbb{E}x_{ij}^{l_1+l_2+l_3}(t)) \tag{4.2}$

where in the first equality, the reason that the degree of each item in the right-hand side is even is that when $r$ is odd, $\mathbb{E}A_{kj}^r = 0$ holds; the conclusion that $\mathbb{E}A_{kj}^r = \frac{a(t)^r}{r+1}$ when $r$ is even is applied to derive the right-hand side of the second equality; and the conclusion that $\mathbb{E}C_{kj}^r = \frac{2^r}{r+1}$ is applied to derive the right-hand side of the last equality. In addition, the expression $\frac{r}{2} - p - q$ is marked as $s$ for symmetry of the expressions above.

According to the Eq. (4.2), $\sigma^r x_{ij}(t+1)$ can be expressed as the sum of several items (corresponding to the five-fold sum $\sum_{p=0}^{\frac{r}{2}} \sum_{q=0}^{\frac{r}{2}-p} \sum_{l_1=0}^{2p} \sum_{l_2=0}^{2q} \sum_{l_3=0}^{2s} *$), where each item can be expressed as the linear combination of $\mathbb{E} x_{ij}^m(t) (m = 0,1, \dots, r)$ (corresponding to the coefficient $\frac{\binom{r}{2p}\binom{r-2p}{2q}}{(2p+1)(2q+1)(2s+1)} \sum_{l_1+l_2+l_3=m} \frac{2^{r-(l_1+l_2+l_3)} \binom{2p}{l} \binom{2q}{l} \binom{2s}{l}}{(2p-l_1+1)(2q-l_2+1)(2s-l_3+1)} p_{1j}^{2p-l_1} p_{1j}^{2q-l_2} p_{1j}^{2s-l_3}$ ). In summary, $\sigma^r x_{ij}(t+1)$ can be expressed as the linear combination of $\mathbb{E} x_{ij}^m(t) (m = 0,1, \dots, r)$.

In addition, $\mathbb{E} x_{ij}^r(t)$ can be calculated as:

$$\mathbb{E} x_{ij}^r(t) = \mathbb{E}\left(\left(x_{ij}(t) - \frac{1}{3}\sum p_j\right) + \frac{1}{3}\sum p_j\right)^r = \sum_{m=0}^{r} \binom{r}{m} \left(\frac{1}{3}\sum p_j\right)^{r-m} \mathbb{E}\left(x_{ij}(t) - \frac{1}{3}\sum p_j\right)^m$$

$$= \sum_{m=0}^{r} \binom{r}{m} \left(\frac{1}{3}\sum p_j\right)^{r-m} \sigma^m x_{ij}(t) \qquad (4.3)$$

According to Eq. (4.3), $\mathbb{E} x_{ij}^r(t)$ can be expressed as the linear combination of $\sigma^m x_{ij}(t) (m = 0,1, \dots, r)$, which means by substituting Eq. (4.3) into Eq. (4.2), $\sigma^r x_{ij}(t+1)$ can be expressed as the linear combination of $\sigma^m x_{ij}(t) (m = 0,1, \dots, r)$ as well. Thus, if the expression of all $m$-order $(m = 0,1, \dots, r)$ central moments of $x_{ij}(t)$ has been derived, we can obtain a difference equation of $\sigma^r x_{ij}(t)$ corresponding to $t$:

$$\sigma^r x_{ij}(t+1) = f_r \sigma^r x_{ij}(t) + g_r(t) \qquad (4.4)$$

where $f_r$ is a coefficient, and $g_r(t)$ denotes the linear combination of $\sigma^m x_{ij}(t) (m = 0,1, \dots, r-1)$.

Eq. (4.4) provides a way of solving the central moment of arbitrary positive integer order of $x_{ij}(t)$ via mathematical induction: from $r = 1$, based on the expressions of $\sigma^m x_{ij}(t) (m = 0,1, \dots, r)$ that have been solved, according to Eq. (4.4), the difference equation that $\sigma^{r+1} x_{ij}(t)$ satisfies corresponding to $t$ is solved. It is worth noting that for any odd positive integer $r$, $\sigma^r x_{ij}(t) = 0$ holds, thus the procedure of solving $\sigma^{r_0} x_{ij}(t)$ ($r_0$ is even) can be further clarified in Fig. 10:

```
1:   σ¹x_ij(t) = 0, E¹x_ij(t) = (1/3)Σp_j
2:   r = 1
3:   while 1
4:       solve σ²ʳx_ij(t) according to Eqs. (4.2) and (4.3)
5:       solve E²ʳx_ij(t) according to Eq. (4.3)
6:       σ^(2r+1) x_ij(t) = 0
7:       solve E^(2r+1) x_ij(t) according to Eq. (4.3)
8:       if 2r + 1 = r_0
9:           break
10:      end if
11:      r = r + 1
12:  end while
13:  return σ^(r_0) x_ij(t)
```

Fig. 10 Pseudo-code of solving $\sigma^{r_0} x_{ij}(t)$ ($r_0$ is even).

As an example, the second order central moment (i.e. variance) of $x_{ij}(t)$ (marked as $\mathbb{D}x_{ij}(t)$) can be solved as:

**Proposition 3.3**

$$\mathbb{D}x_{ij}(t+1) = \frac{a(t)^2}{9}\left(\mathbb{D}x_{ij}(t) + \frac{4}{9}\sum p_j^2 - \frac{1}{9}\left(\sum p_j\right)^2\right) \quad (4.5)$$

In the rest part of this section, the characteristics of $\mathbb{D}x_{ij}(t)$ will be further discussed, which could serve as tools for the stability analysis of the GWO to be carried out in the next section. For simplification, $\mathbb{D}x_{ij}(t)$ is abbreviated as $D_t$, $\frac{a(t)^2}{9}$ is abbreviated as $b_t$ (i.e. $b_t = \frac{4}{9}\left(1-\frac{t}{T}\right)^2$), and expression $\frac{4}{9}\sum p_j^2 - \frac{1}{9}\left(\sum p_j\right)^2$ is marked as $p_0$ (obviously $p_0 = \frac{1}{9}\sum p_j^2 + \frac{1}{9}\sum(p_i-p_j)^2 \geq 0$ holds), then Eq. (4.5) can be rewritten as:

$$D_{t+1} = b_t(D_t + p_0) \quad (4.6)$$

The Eq. (4.6) can be seen as a first-order dynamic system. We firstly suppose that $b_t$ is constant (independent of $t$), since $\forall t = 1, \ldots, T, 0 \leq b_t < 1$ holds, this dynamic system holds a unique attractor $D_0 = \frac{b_t p_0}{1-b_t}$, which means if $b_t$ is constant, $\lim_{t\to\infty} D_t = D_0$. In addition, it is worth noting that this assumption is approximately satisfied when the value of $T$ is large enough (in which case $|b_{t+1} - b_t| = \frac{4|2T-2t-1|}{9T^2}$ tends to 0 as $T$ increases), thus it could be assumed that when $T \to \infty$, $\lim_{t\to\infty} D_t = D_0$ holds. The analysis above can be summarized as follows:

**Corollary 3.1**

When $T \to \infty$, $\left|D_t - \frac{b_t p_0}{1-b_t}\right| \leq \left(\frac{4}{9}\right)^{t-1}\left|D_1 - \frac{b_1 p_0}{1-b_1}\right|$ holds.

According to the corollary 3.1, an upper bound of the distance between $D_t$ and $D_0 = \frac{b_t p_0}{1-b_t}$ can be obtained when $T \to \infty$, by which we can deduce that when $T \to \infty$, $\lim_{t\to\infty} D_t = D_0$ holds.

The conclusions above can be verified via MC simulations. Here, we set $T = 50$, and the values of $p_{kj}$ are randomly chosen. For each group of $p_{kj}$ chosen, the GWO under the stagnation assumption is run for $10^6$ trials independently (total iteration $T = 50$). If the value of $x_{ij}$ obtained in the $l$th trial with iteration $t$ is marked as $x_{ij}(l,t)$, and $\forall l, x_{ij}(l,1)$ is chosen as a random variable following $U[-4,4]$, then $D_t$ can be approximated as the variance of $x_{ij}(l,t), l = 1, \ldots, 10^6$. The $D_t$ obtained from MC simulations and $D_0 = \frac{b_t p_0}{1-b_t}$ are compared as shown in Fig. 11.

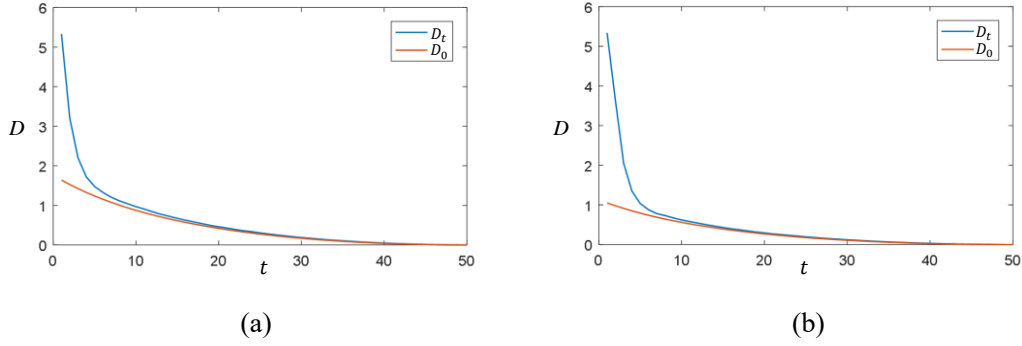

(a) (b)

Fig. 11 $D_t$ and $D_0$, (a): $p_0 = 2.20$, (b): $p_0 = 1.41$

It can be seen from Fig. 11 that with the increasing of $t$, the value of $D_t$ tends to $D_0$. In addition, for each $T$, the MC simulation is carried out, and $|D_T - D_0|$ is calculated and shown in Table II, from which we can see that the value of $|D_T - D_0|$ is close to 0, and decreases with the increasing of $T$, which verifies the conclusion that $\lim_{T \to \infty} D_T = D_0$.

Table II Relationship between $|D_T - D_0|$ and $T$

| $T$ | 20 | 30 | 50 | 100 | 200 | 500 |
|---|---|---|---|---|---|---|
| $|D_T - D_0|$ | $2.3 \times 10^{-3}$ | $1.7 \times 10^{-3}$ | $4.6 \times 10^{-4}$ | $3.7 \times 10^{-5}$ | $7.1 \times 10^{-6}$ | $2.2 \times 10^{-6}$ |

By calculating and analyzing the moments of $x_{ij}(t)$, the characteristics of the updating mechanism of the GWO can be further understood. In the part II of the paper, we will prove the global convergence of the GWO based on the results of the moments of $x_{ij}(t)$ calculated above. In addition, as another example, in the next section, the stability analysis of the GWO under the stagnation assumption will be carried out.

**4.3 Stability analysis of GWO under the stagnation assumption**

As an application of the results of $\mathbb{E}x_{ij}(t)$ and $\mathbb{D}x_{ij}(t)$ derived in the section 4.2, the stability analysis of GWO under the stagnation assumption is proved and discussed. As introduced in section 1, as an important branch of the theoretical analysis of meta-heuristic algorithms, the stability analysis has been carried out deeply, aiming at studying the dynamics of search agents, and moreover, on what kind of conditions do the search agents converge to some certain position [59].

Taking the PSO as an example, early analysis often ignores the randomness and sets the random coefficients as constants (whose values are often set as their expectations) [60-61], by which the updating mechanism of PSO could be simplified into a deterministic difference equation, thus the stability analysis could be performed in the deterministic dynamic systems by using methods such as Lyapunov stability analysis [62]. Obviously, the assumption that all the random coefficients are constants is too strong to fit the actual situation well. Subsequently, the stability analysis performed on the PSO with randomness consideration is widely studied. When randomness presents, a popular substitute called order-1 stability is defined as [46]:

**Definition 1**: A meta-heuristic algorithm holds order-1 stability if for any search agent $i$, $\lim_{t \to \infty} \mathbb{E}x_{ij}(t)$ exists.

Obviously, in the concept of order-1 stability, the deterministic model is established by replacing the random variable $x_{ij}(t)$ with $\mathbb{E}x_{ij}(t)$.

Despite the fruitful research achievements of the order-1 stability of the PSO [62-64], it has been pointed out that such order-1 stability is not enough to ensure convergence, and the following definition of order-2 stability should be taken into account [46,65]:

**Definition 2**: A meta-heuristic algorithm holds order-2 stability if for any search agent $i$, $\lim_{t \to \infty} \mathbb{D}x_{ij}(t) = 0$.

It is worth noting that till now, for simplicity, the stability analysis of the PSO is mostly carried out under the stagnation assumption, in which the global best (gbest) of the population and the personal best (pbest) of each search agent are set to be constant (as mentioned in section 3.1).

Based on the analysis above, as a direct corollary of the proposition 3.2, the stability analysis of GWO under the stagnation assumption is performed as follows:

**Proposition 3.4**

The GWO holds order-1 stability under the stagnation assumption.

and

**Proposition 3.5**

The GWO holds order-2 stability under the stagnation assumption.

In the last part of this section, the two propositions above are verified by numerical simulations. The simulation method is identical with that of obtaining $D_t$ and $D_0$ in Fig. 11 in section 4.2, i.e., $T = 50$, four groups of $p_{kj}$ are randomly chosen, and for each group of $p_{kj}$ chosen, the GWO under stagnation assumption is run for $10^6$ trials independently. Here, we mark the value of $x_{ij}$ obtained in the $l$th trial with iteration $t$ as $x_{ij}(l,t)$ (in each trial the value of $x_{ij}(l,1)$ is randomly chosen within $[-4,4]$, corresponding to the optimization problem with feasible region $[-4,4]$), then $\mathbb{E}x_{ij}(t)$ could be approximated as the mean of $x_{ij}(l,t)$ (i.e. $10^{-6}\sum_{l=1}^{10^6} x_{ij}(l,t)$), and $\mathbb{D}x_{ij}(t)$ could be approximated as the variance of $x_{ij}(l,t)$ (i.e. $\frac{1}{10^6-1}\sum_{l=1}^{10^6}(x_{ij}(l,t) - 10^{-6}\sum_{m=1}^{10^6} x_{ij}(m,t))^2$. It is worth noting that after the stagnation assumption, whatever the value of $f(\mathbf{x_i}(t))$ is, the position of the best three search agents will not be updated, which means the information of objective function is no longer required during simulations. The simulation results of $\mathbb{E}x_{ij}(t)$ and $\mathbb{D}x_{ij}(t)$ are shown in Fig. 12 and Fig. 13 (in which the Y-axis is set as $\ln \mathbb{D}x_{ij}(t)$ to better exhibit the decreasing speed of $\mathbb{D}x_{ij}(t)$), respectively.

It can be seen from Fig. 12 that $\mathbb{E}x_{ij}(t) = \begin{cases} 0, t = 1 \\ \frac{1}{3}\sum p_j, \text{otherwise} \end{cases}$, which is a direct inference of $x_{ij}(t) \sim U[-4,4]$ and proposition 3.2. In addition, it can be seen from Fig. 13 that when $t \to T$, the value of $\ln \mathbb{D}x_{ij}(t)$ is negative and decreases sharply, thus the conclusion $\lim_{t \to \infty} \mathbb{D}x_{ij}(t) = 0$ can be verified.

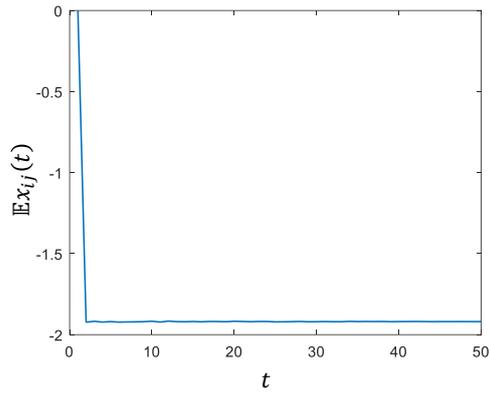

(a)

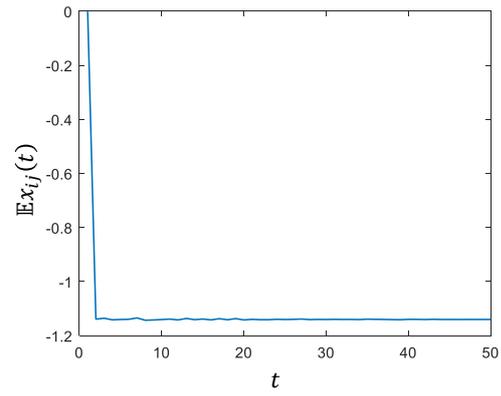

(b)

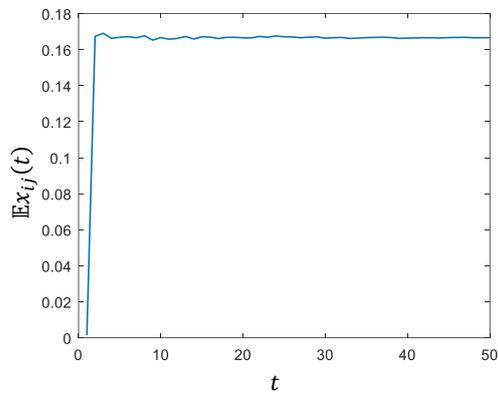

(c)

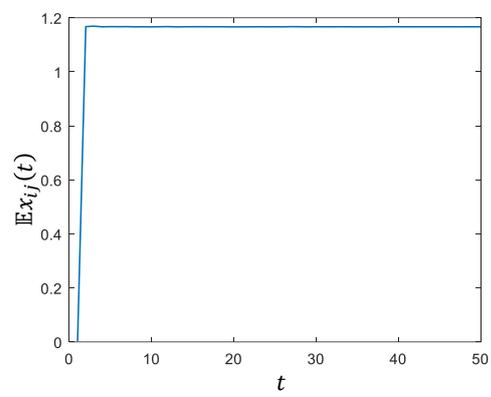

(d)

Fig. 12 Simulation results of $\mathbb{E}x_{ij}(t)$, (a): $\frac{1}{3}\sum p_j = -1.92$, (b): $\frac{1}{3}\sum p_j = -1.14$, (c): $\frac{1}{3}\sum p_j = 0.17$, (d): $\frac{1}{3}\sum p_j = -0.28$

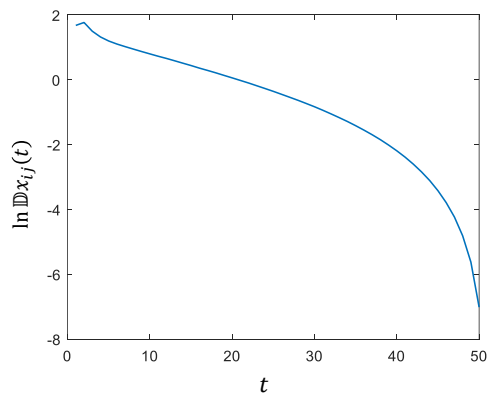

(a)

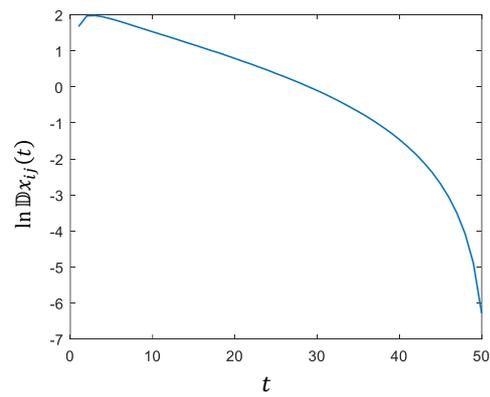

(b)

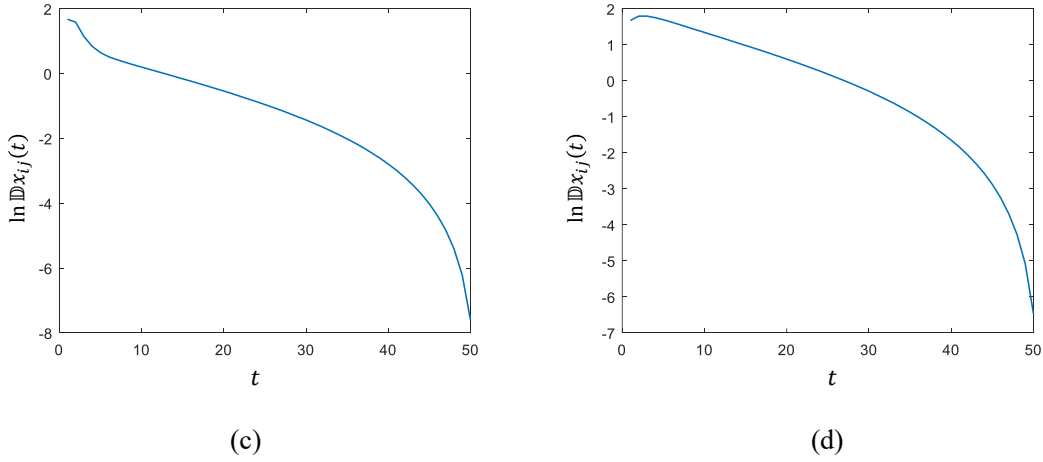

Fig. 13 Simulation results of $\ln \mathbb{D}x_{ij}(t)$, (a): $\frac{1}{3}\sum p_j = 0.84$, (b): $\frac{1}{3}\sum p_j = 2.14$, (c): $\frac{1}{3}\sum p_j = -1.89$, (d): $\frac{1}{3}\sum p_j = 2.20$

## 5. Conclusion

In the part I of the paper, we establish a theoretical framework of the GWO and discuss several promising theoretical findings, containing the sampling distribution, the order-1 and order-2 stability, and the global convergence analysis. According to the stagnation assumption, that is, the positions of the best three search agents $\mathbf{p_k}(t)(k = 1,2,3)$ remain unchanged during iterations, the characteristics of the sampling distribution of the new solution $\mathbf{x_i}(t + 1)$ and the probabilistic stability of the GWO are elaborately addressed.

For simplification, the characteristics of the sampling distribution of the new solution $\mathbf{x_i}(t + 1)$ is firstly discussed under the assumption that the original solution before updating $\mathbf{x_i}(t)$ is constant, in which case $\mathbf{x_i}(t + 1)$ becomes a function of random variables $A_{kj}, C_{kj}(k = 1,2,3)$. Although the exact analytical expression of $x_{ij}(t + 1)$'s PDF (marked as $h_{ij}(u_j)$) is hard to derive, it can be proved that $h_{ij}(u_j)$ is a continuous single-peak function symmetrical about $u_j = \frac{1}{3}\sum p_j$ with a bounded domain. In addition, the differentiability of $h_{ij}(u_j)$ is carefully studied.

Regarding more realistic condition that $\mathbf{x_i}(t)$ is no longer regarded as a constant, the conclusion that the $h_{ij}(u_j)$ is a single-peak function symmetrical about $u_j = \frac{1}{3}\sum p_j$ still holds. In addition, the central moments of any positive integer order of $x_{ij}(t)$ is analyzed, and the order-1 and order-2 stability of the GWO under stagnation assumption is also proved. All the above analysis results are verified by numerical simulations.

The corresponding results obtained in the part I will serve as the basis for the global convergence analysis of the GWO in the part II. At the same time, we also believe that the theoretical findings in the part I of the paper can be used as a fundamental tool to some extent to help researchers better understand the searching mechanism and particle trajectory of the GWO, and future research on the application, novel variants, and more in-depth theoretical findings of the GWO can be carried out based on our works.

## Appendix

**Proof of proposition 1.2**

Mark $g_{kj}(u_j)$ as $g_{kj}(u_j, x_{ij}, p_{kj})$, then according to Eq. (2.2), $g_{kj}(u_j, x_{ij}, p_{kj}) = g_{kj}(u_j, -x_{ij}, -p_{kj}) = g_{kj}(u_j, 2p_{kj} - x_{ij}, p_{kj}) = g_{kj}(2p_{kj} - u_j, x_{ij}, p_{kj})$, thus we only need to consider the conditions $p_{kj} \geq 0$, $x_{ij} \geq p_{kj}$, and $u_j \geq p_{kj}$.

Let $D = \{(A_{kj}, C_{kj}) | p_{kj} + A_{kj} | C_{kj} p_{kj} - x_{ij}| \leq u_j, -a(t) \leq A_{kj} \leq a(t), 0 \leq C_{kj} \leq 2\}$, then according to Eq. (2.2), the cumulative distribution function (CDF) of $x'_{kj}$ can be derived as

$$G_{kj}(u_j) = P\{x'_{kj} \leq u_j\} = \iint_D \frac{1}{4a(t)} dA_i dC_i.$$

When $x_{ij} > 2p_{kj}$ and $u_j > a(t)(x_{ij} - 2p_{kj})$, region $D$ is shown in Fig. 14 (a) (marked by grey shading, similarly hereinafter). In this case, we have

$$G_{ki}(u_j) = \frac{1}{4a(t)}\left[2a(t) + \frac{2u_j}{x_{ij}} + \int_{\frac{u_j}{x_{ij}}}^{a(t)} \frac{u_j}{p_{kj}} \ln v \, dv - \left(\frac{x_{ij}}{p_{kj}} - 2\right)\left(a(t) - \frac{u_j}{x_{ij}}\right)\right]$$

$$= \frac{1}{2} + \frac{1}{4a(t)}\left[\frac{u_j}{p_{kj}} + \frac{u_j}{p_{kj}} \ln \frac{a(t)x_{ij}}{u_j} - a(t)\left(\frac{x_{ij}}{p_{kj}} - 2\right)\right]$$

When $x_{ij} > 2p_{kj}$ and $0 \leq u_j \leq a(t)(x_{ij} - 2p_{kj})$, region $D$ is shown in Fig. 14 (b). In this case, we have

$$G_{kj}(u_j) = \frac{1}{2} + \frac{u_j}{4a(t)p_{kj}} \ln \frac{x_{ij}}{x_{ij} - 2p_{kj}}$$

When $p_{kj} \leq x_{ij} \leq 2p_{kj}$ and $u_j > a(t)(2p_{kj} - x_{ij})$, region $D$ is shown in Fig. 14 (c). In this case, we have

$$G_{kj}(u_j) = \frac{1}{2} + \frac{1}{4a(t)}\left[\frac{u_j}{p_{kj}} + \frac{u_j}{p_{kj}} \ln \frac{a(t)x_{ij}}{u_j} + a(t)\left(2 - \frac{x_{ij}}{p_{kj}}\right)\right]$$

When $p_{kj} \leq x_i \leq 2p_{kj}$ and $0 \leq u_j \leq a(t)(2p_{kj} - x_{ij})$, region $D$ is shown in Fig. 14 (d). In this case, we have

$$G_{kj}(u_j) = \frac{1}{2} + \frac{1}{2a(t)}\left[\frac{u_j}{p_{kj}} + \frac{u_j}{p_{kj}} \ln \frac{a(t)\sqrt{x_{ij}(2p_{kj} - x_{ij})}}{u_j}\right]$$

Thus

$$g_{kj}(u_j) = \frac{dG_{kj}(u_j)}{du_j} =$$

$$\begin{cases} \frac{1}{4a(t)p_{kj}} \ln \frac{a(t)x_{ij}}{u_j}, & x_{ij} > 2p_{kj}, u_j > a(t)(x_{ij} - 2p_{kj}) \\ \frac{1}{4a(t)p_{kj}} \ln \frac{x_{kj}}{x_{kj} - 2p_{kj}}, & x_{ij} > 2p_{kj}, 0 \leq u_j \leq a(t)(x_{ij} - 2p_{kj}) \\ \frac{1}{4a(t)p_{kj}} \ln \frac{a(t)x_{ij}}{u_j}, & p_{kj} \leq x_{ij} \leq 2p_{kj}, u_j > a(t)(2p_{kj} - x_{ij}) \\ \frac{1}{2a(t)p_{kj}} \ln \frac{a(t)\sqrt{x_{ij}(2p_{kj} - x_{ij})}}{u_j}, & p_{kj} \leq x_i \leq 2p_{kj}, 0 \leq u_j \leq a(t)(2p_{kj} - x_{ij}) \end{cases}$$

Expanding the equations above to other cases where $x_{ij}, p_{kj}, u_j \in \mathbb{R}$, and after a simple arrangement with absolute value sign, we finally get Eqs. (3.4) - (3.6). □

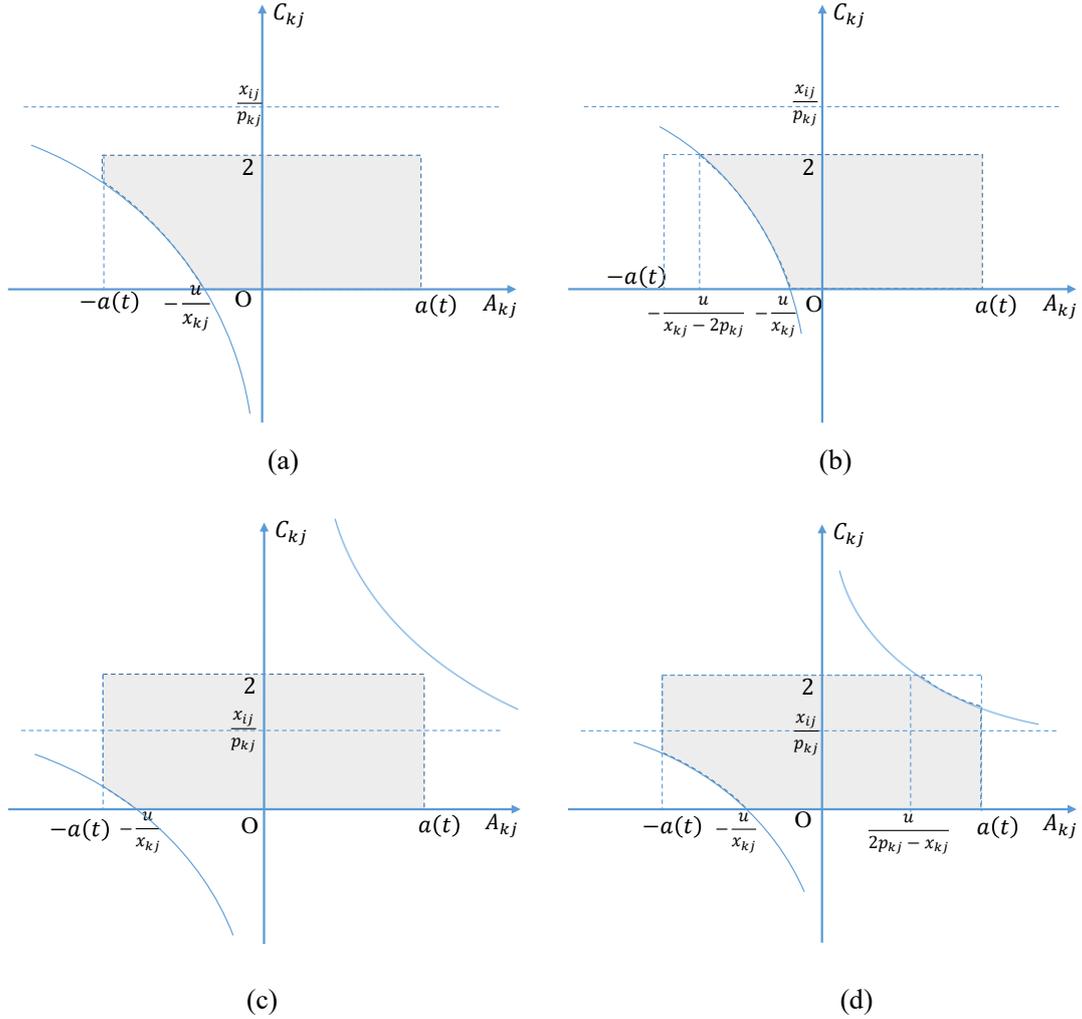

Fig. 14 Region $D$ (marked by grey shading), (a): $x_{ij} > 2p_{kj}$ and $u_j > a(t)(x_{ij} - 2p_{kj})$; (b): $x_{ij} > 2p_{kj}$ and $0 \leq u_j \leq a(t)(x_{ij} - 2p_{kj})$; (c): $p_{kj} \leq x_{ij} \leq 2p_{kj}$ and $u_j > a(t)(2p_{kj} - x_{ij})$; (d): $p_{kj} \leq x_{ij} \leq 2p_{kj}$ and $0 \leq u_j \leq a(t)(2p_{kj} - x_{ij})$.

**Proof of corollary 2.1**

According to the characteristics of convolution operator, the domain of $h_{ij}(u_j)$ is evident.

It is well known that the convolution operator holds evident smoothing effect, for instance, the convolution of a continuous function and an arbitrary function must be a continuous function as well [57], i.e., the curve of the arbitrary function selected (which might not be a continuous function) is smoothed. It is worth noting that when $m_{kj} > 0$, $g_{kj}(u_j)$ is continuous within $\mathbb{R}$, therefore we only need to consider the case that $m_{kj} \leq 0 (k = 1,2,3)$. Since $g_{kj}(u_j)$ is continuous within $\mathbb{R}\setminus\{p_{kj}\}$ and $\lim_{u \to p_{kj}} g_{kj}(u_j) = +\infty$ holds (please see proposition 1.2), thus the continuity of $(g_{1j} * g_{2j})(u_j)$ can be proved only $(g_{1j} * g_{2j})(p_{1j} + p_{2j}) < +\infty$ holds.

When $m_{kj} \leq 0 (k = 1,2,3)$, just let $|m_{1j}| \leq |m_{2j}|$, then $(g_{1j} * g_{2j})(p_{1j} + p_{2j})$ can be calculated as:

$(g_{1j} * g_{2j})(p_{1j} + p_{2j})$
$= 2\int_{p_{1j}-n_{1j}}^{p_{1j}+m_{1j}} g_{1j}(u_j)g_{2j}(p_{1j} + p_{2j} - u_j)du + 2\int_{p_{1j}+m_{1j}}^{p_{1j}} g_{1j}(u_j)g_{2j}(p_{1j} + p_{2j} - u_j)du_j$

$$< 2(n_{1j} + m_{1j}) \cdot \frac{\ln\frac{\sqrt{-m_{1j}n_{1j}}}{-m_{1j}}}{n_{1j}-m_{1j}} \cdot \frac{\ln\frac{\sqrt{-m_{2j}n_{2j}}}{-m_{2j}}}{n_{2j}-m_{2j}} +$$

$$\frac{2}{(n_{1j}-m_{1j})(n_{2j}-m_{2j})}\int_0^{-m_{1j}} \ln\frac{\sqrt{-m_{1j}n_{1j}}}{u_j} \ln\frac{\sqrt{-m_{2j}n_{2j}}}{u_j} du_j$$

$$= A + k\left\{\ln P \ln Q - (\ln P + \ln Q)(u_j \ln u_j - u_j) + \left(u_j(\ln u_j)^2 - 2u_j \ln u_j + 2u_j\right)\right\}\Big|_0^{-m_{1j}}$$

$$= A + k\left\{-m_{1j}\ln P \ln Q - m_{1j}(\ln P + \ln Q)(1 - \ln(-m_{1j})) - m_{1j}\left((\ln(-m_{1j}))^2 - \right.\right.$$

$$\left.\left. 2\ln(-m_{1j}) + 2\right)\right\} < +\infty$$

where $A = 2(n_{1j} + m_{1j}) \cdot \frac{\ln\frac{\sqrt{-m_{1j}n_{1j}}}{-m_{1j}}}{n_{1j}-m_{1j}} \cdot \frac{\ln\frac{\sqrt{-m_{2j}n_{2j}}}{-m_{2j}}}{n_{2j}-m_{2j}}$, $k = \frac{2}{(n_{1j}-m_{1j})(n_{2j}-m_{2j})}$, $P = \sqrt{-m_{1j}n_{1j}}$, $Q = \sqrt{-m_{2j}n_{2j}}$ are coefficients.

Thus the continuity of $(g_{1j} * g_{2j})(u_j)$ is proved, i.e., the continuity of $h_{ij}(u_j)$ is proved.□

**Proof of corollary 2.2**

The first two conclusions are evident due to the characteristics of convolution operator. As for the latter part, we first prove that $(g_{1j} * g_{2j})(u_j)$ is non-decreasing within $(-\infty, p_{1j} + p_{2j}]$.

For the sake of simplification, let $g_{10}(u_j) = g_{1j}(u_j + p_{1j})$, $g_{20}(u_j) = g_{2j}(u_j + p_{2j})$, then according to the corollary 1.1, both $g_{10}(u_j), g_{20}(u_j)$ are non-decreasing on $\mathbb{R}^-$, and we only need to prove that $h_0(u_j) = (g_{10} * g_{20})(u_j)$ is non-decreasing on $\mathbb{R}^-$. Let $u_j < u_j + \Delta u_j < 0$, we have

$$h_0(u_j + \Delta u_j) - h_0(u_j)$$
$$= \int_{-\infty}^{+\infty} g_{10}(v)\left(g_{20}(u_j + \Delta u_j - v) - g_{20}(u_j - v)\right) dv$$

Let $c(v) = g_{20}(u_j + \Delta u_j - v) - g_{20}(u_j - v)$, then obviously $c(2u_j + \Delta u_j - v) = -c(v)$, and when $v < u_j + \frac{1}{2}\Delta u_j$, $c(v) \leq 0$. Thus

$$h_0(u_j + \Delta u_j) - h_0(u_j)$$
$$= \int_{-\infty}^{u_j + \frac{1}{2}\Delta u_j} g_{10}(v)c(v)\, dv + \int_{u_j+\frac{1}{2}\Delta u_j}^{+\infty} g_{10}(v)c(v)\, dv$$

$$= \int_{-\infty}^{u_j + \frac{1}{2}\Delta u_j} \left(g_{10}(v) - g_{10}(2u_j + \Delta u_j - v)\right)c(v)\, dv \geq 0$$

i.e. $h_0(u_j)$ is non-decreasing on $\mathbb{R}^-$, and consequently $(g_{1j} * g_{2j})(u_j)$ is non-decreasing within $(-\infty, p_{1j} + p_{2j})$. Similarly, $h_{ij}(u_j)$ is non-decreasing within $\left(-\infty, \frac{1}{3}\sum_{k=1}^3 p_{kj}\right)$.□

**Proof of corollary 2.3**

Just assume that $m_{kj} > 0$, in which case $g_{3j}(u_j)$ is continuous and bounded within $\mathbb{R}$. It can be deduced from the corollary 2.1 that whatever the signs of $m_{kj}(k = 1,2)$ are, $(g_{1j} * g_{2j})(u_j)$ is always a continuous and bounded function within $\mathbb{R}$. Mark $(g_{1j} * g_{2j})(u_j)$ as

$d(u_j)$, then we only need to prove that $\frac{d}{du}(d * g_{3j})(u_j)$ exists and is continuous within $\mathbb{R}$.

For the convenience of analysis, all subsection region endpoints of $d(u_j)$ are set as (from left to right) $[s_0, s_1, ..., s_{q_1}]$, and all subsection region endpoints of $g_{3j}(u_j)$ are set as (from left to right) $[t_0, t_1, ..., t_{q_2}]$. By introducing unit step function $\varepsilon(x) = \begin{cases} 0, x \leq 0 \\ 1, x > 0 \end{cases}$, the expressions of $g_{1j}(u_j)$ and $g_{2j}(u_j)$ can be rewritten as:

$$d(u_j) = \sum_{l_1=1}^{q_1}[\varepsilon(u_j - s_{l_1-1}) - \varepsilon(u_j - s_{l_1})]v_{l_1}(u_j)$$
$$g_{3j}(u_j) = \sum_{l_2=1}^{q_2}[\varepsilon(u_j - t_{l_2-1}) - \varepsilon(u_j - t_{l_2})]w_{l_2}(u_j)$$

where $v_{l_1}(u_j), w_{l_2}(u_j)$ are continuous functions defined within $\mathbb{R}$, $\forall l_1 = 1, ..., q_1, \forall u_j \in [s_{l_1-1}, s_{l_1}], v_{l_1}(u_j) = g_{1j}(u_j)$; $\forall l_2 = 1, ..., q_2, \forall u_j \in [t_{l_2-1}, t_{l_2}], w_{l_2}(u_j) = g_{2j}(u_j)$, and $\forall l_1 = 1, ..., q_1 - 1, v_l(s_l) = v_{l+1}(s_l)$ holds.

The expression of $(d * g_{3j})(u_j)$ can be expanded as:

$(d * g_{3j})(u_j)$

$= \int_{-\infty}^{+\infty}\{\sum_{l_1=1}^{q_1}[\varepsilon(z - s_{l_1-1}) - \varepsilon(z - s_{l_1})]l_1(z)\}\{\sum_{l_2=1}^{q_2}[\varepsilon(u_j - z - t_{l_2-1}) - \varepsilon(u_j - z - t_{l_2})]w_{l_2}(u_j - z)\}dz$

$= \sum_{l_1=1}^{q_1}\sum_{l_2=1}^{q_2}\int_{-\infty}^{+\infty}[\varepsilon(z - s_{l_1-1}) - \varepsilon(z - s_{l_1})]\varepsilon(u_j - z - t_{l_2-1})v_{l_1}(z)w_{l_2}(u_j - z)dv -$

$\int_{-\infty}^{+\infty}[\varepsilon(z - s_{l_1-1}) - \varepsilon(z - s_{l_1})]\varepsilon(u_j - z - t_{l_2})v_{l_1}(z)w_{l_2}(u_j - z)dz$

Thus

$\frac{d}{du_j}(d * g_{3j})(u_j)$

$= \sum_{l_1=1}^{q_1}\sum_{l_2=1}^{q_2}\int_{-\infty}^{+\infty}[\varepsilon(z - s_{l_1-1}) - \varepsilon(z - s_{l_1})]\varepsilon(u_j - z - t_{l_2-1})v_{l_1}(z)w'_{l_2}(u_j - z)dz +$

$[\varepsilon(u_j - t_{l_2-1} - s_{l_1-1}) - \varepsilon(u_j - t_{l_2-1} - s_{l_1})]v_{l_1}(u_j - t_{l_2-1})w_{l_2}(t_{l_2-1}) - \int_{-\infty}^{+\infty}[\varepsilon(z - s_{l_1-1}) - \varepsilon(z - s_{l_1})]\varepsilon(u_j - z - t_{l_2})v_{l_1}(z)w'_{l_2}(u_j - z)dz - [\varepsilon(u_j - t_{l_2} - s_{l_1-1}) - \varepsilon(u_j - t_{l_2} - s_{l_1})]v_{l_1}(u_j - t_{l_2})w_{l_2}(t_{l_2})$

$= (d * g'_{3j})(u_j) + \sum_{l_1=1}^{q_1}\sum_{l_2=1}^{q_2}[\varepsilon(u_j - t_{l_2-1} - s_{l_1-1}) - \varepsilon(u_j - t_{l_2-1} - s_{l_1})]v_{l_1}(u_j - t_{l_2-1})w_{l_2}(t_{l_2-1}) - [\varepsilon(u_j - t_{l_2} - s_{l_1-1}) - \varepsilon(u_j - t_{l_2} - s_{l_1})]v_{l_1}(u_j - t_{l_2})w_{l_2}(t_{l_2})$

where $g'_{3j}(u_j)$ denotes the piecewise function in which each subsection $w_{l_2}(u_j)$ is replaced by $w'_{l_2}(u_j) = \frac{dw'_{l_2}(u_j)}{du_j}$, rather than the derivative of $g_{3j}(u_j)$ (in fact $g_{3j}(u_j)$ is non-derivable on its subsection region endpoints).

Since $d(u_j)$ is continuous and bounded within $\mathbb{R}$, the first item of $\frac{d}{du_j}(d * g_{3j})(u_j)$ (i.e. $(d * g'_{3j})(u_j)$) is continuous within $\mathbb{R}$ as well. Since $g_{3j}(u_j)$ is bounded within $\mathbb{R}$ and $\forall l_1 = 1, ..., q_1 - 1, v_l(s_l) = v_{l+1}(s_l)$ holds, the second item of $\frac{d}{du_j}(d * g_{3j})(u_j)$ is continuous within

$\mathbb{R}$ as well. In summary, when there exists at least one positive integer among $m_{kj}(k = 1,2,3)$, $h_{ij}(u_j)$ is differentiable within $\mathbb{R}$.

In addition, since $h_{ij}\left(u_j + \frac{1}{3}\sum p_j\right)$ is an even function, $h'_{ij}\left(u_j + \frac{1}{3}\sum p_j\right)$ becomes an odd function. It is worth noting that $h'_{ij}\left(u_j + \frac{1}{3}\sum p_j\right)$ is continuous within $\mathbb{R}$, thus we have $h'_{ij}\left(\frac{1}{3}\sum p_j\right) = 0$. □

**Proof of proposition 3.1**

Obviously, $\forall t \geq 2$, the domain of $h_{ij}(u_j, t)$ is symmetrical about $\frac{1}{3}\sum p_j$. Mark the right endpoint of the domain of $h_{ij}(u_j, t)$ as $u(t)$, if $x_{ij}(1)$ is fixed, then according to corollary 2.1, we have:

$$u(2) - \frac{1}{3}\sum p_j = \frac{1}{3}\sum n_j = \frac{a(1)}{3}\left[\sum_{k=1}^{3}|p_k| + \sum_{k=1}^{3}|x_{ij}(1) - p_k|\right] \geq \frac{a(1)}{3}\sum_{k=1}^{3}|x_{ij}(1) - p_k| \geq a(1)\left|x_{ij}(1) - \frac{1}{3}\sum p_j\right|.$$

Taking the fact that $x_{ij}(1)$ is a random variable into consideration, we have $u(2) - \frac{1}{3}\sum p_j \geq a(1)\max_{x_{ij}(1)}\left|x_{ij}(1) - \frac{1}{3}\sum p_j\right| = a(1)d_1$.

By repeating the derivations above, we can get:

$$u(t) - \frac{1}{3}\sum p_j \geq a(t-1)\max_{x_{ij}(t-1)}\left|x_{ij}(t-1) - \frac{1}{3}\sum p_j\right| = a(t-1)\left(u(t-1) - \frac{1}{3}\sum p_j\right) \geq a(t-1)a(t-2)\left(u(t-2) - \frac{1}{3}\sum p_j\right) \geq \cdots \geq d_1 \prod_{\tau=1}^{t-1} a(\tau)$$

which means the domain of $h_{ij}(u_j, t)$ contains $\left(\frac{1}{3}\sum p_j - d_1 \prod_{\tau=1}^{t-1} a(\tau), \frac{1}{3}\sum p_j + d_1 \prod_{\tau=1}^{t-1} a(\tau)\right)$. □

**Proof of proposition 3.2**

Eq. (2.2) shows that $x'_{kj}(t) = p_{kj} + A_{kj}|C_{kj}p_{kj} - x_{ij}(t)|$, where $A_{kj} \sim U[-a, a]$, i.e. the distribution that $A_{kj}$ follows is symmetrical about $u_j = 0$. Since $A_{kj}$ is independent of $|C_{kj}p_{kj} - x_{ij}(t)|$, whatever distribution $x_{ij}(t)$ follows, the distribution that $x'_{kj}(t)$ follows must be symmetrical about $\frac{1}{3}\sum p_j$, which means the curve of $h_{ij}(u_j, t)$ is symmetrical about $u_j = \frac{1}{3}\sum p_j$. Next, we prove that $h_{ij}(u_j, t)$ is non-decreasing within $\left(-\infty, \frac{1}{3}\sum p_j\right)$ via mathematical induction.

(1) when $t = 1$, according to the procedure of the GWO, $x_{ij}(1)$ follows a uniform distribution within the feasible region, which makes the conclusion evident.

(2) if the conclusion is correct when $t = t_0$, i.e. $h_{ij}(u_j, t_0)$ is non-decreasing within $\left(-\infty, \frac{1}{3}\sum p_j\right)$, then according to the corollary 2.2, as the sum of two random variables $-x_{ij}(t_0)$ and $C_{kj}p_{kj}$ that both follow single-peak symmetrical distributions, $C_{kj}p_{kj} - x_{ij}(t_0)$ follows a single-peak symmetrical distribution as well. Since $A_{kj}$ is independent of $C_{kj}p_{kj} - x_{ij}(t_0)$, and the distribution that $A_{kj}$ follows is symmetrical about $u_j = 0$, it can be seen that $A_{kj}|C_{kj}p_{kj} - x_{ij}(t_0)|$ and $A_{kj}\left(C_{kj}p_{kj} - x_{ij}(t_0)\right)$ follow the same distribution, in which case we only need to prove that the PDF of $A_{kj}\left(C_{kj}p_{kj} - x_{ij}(t_0)\right)$ is non-increasing within $(0, +\infty)$.

Let the PDF of $A_{kj}\left(C_{kj}p_{kj} - x_{ij}(t_0)\right)$ be $\varphi(u_j)$, by symmetry, just suppose that the symmetrical axis of the curve of $\varphi(u_j)$ is on the left of Y-axis, which means $\varphi(u_j)$ is non-increasing within $(0, +\infty)$; mark the PDF of $A_{kj}$ as $\theta(u_j) = \frac{1}{2a}\left(\varepsilon(u_j + a) - \varepsilon(u_j - a)\right)$, the PDF and CDF of $A_{kj}\left(C_{kj}p_{kj} - x_{ij}(t_0)\right)$ as $\delta(u_j)$ and $\Delta(u_j)$ respectively, then the expression of $\Delta(u_j)$ can be derived as:

$$\Delta(u_j) = \iint_{vw \leq u_j} \varphi(v)\theta(w)dvdw = \frac{1}{2} + 2\int_0^a \frac{1}{2a}dw \int_0^{\frac{u_j}{w}} \varphi(v)dv.$$

Thus $\delta(u_j) = \frac{d\Delta(u_j)}{du_j} = \frac{1}{a}\int_0^a w\varphi\left(\frac{u_j}{w}\right)dw$ holds.

Let $\Delta u_j > 0$, then $\delta(u_j + \Delta u_j) - \delta(u_j) = \frac{1}{a}\int_0^a w\left[\varphi\left(\frac{u_j + \Delta u_j}{w}\right) - \varphi\left(\frac{u_j}{w}\right)\right]dw$. Since $\varphi(u_j)$ is non-increasing within $(0, +\infty)$, $\forall w > 0, u_j > 0, \varphi\left(\frac{u_j + \Delta u_j}{w}\right) - \varphi\left(\frac{u_j}{w}\right) \leq 0$ holds, thus $\delta(u_j + \Delta u_j) - \delta(u_j) \leq 0$ holds. Owing to the arbitrariness of positive numbers $u_j, \Delta u_j$, we can deduce that $\delta(u)$ is non-increasing within $(0, +\infty)$.

Finally, according to the corollary 2.2, we can prove that the PDF of $x_{ij}(t_0 + 1) = \frac{1}{3}\sum p_j + \frac{1}{3}\sum_{k=1}^3 A_{kj}|C_{kj}p_{kj} - x_{ij}(t_0)|$ is non-decreasing within $\left(-\infty, \frac{1}{3}\sum p_j\right)$, which means $\forall t, h_{ij}(u_j, t)$ is non-decreasing within $\left(-\infty, \frac{1}{3}\sum p_j\right)$. □

**Proof of proposition 3.3**

According to Eq. (4.2):

$$\mathbb{D}x_{ij}(t+1) = \frac{1}{9}\mathbb{E}\left(\sum_{k=1}^3 A_{kj}|C_{kj}p_{kj} - x_{ij}(t)|\right)^2 = \frac{1}{9}\sum_{k=1}^3 \mathbb{E}A_{kj}^2 \mathbb{E}\left(C_{kj}p_{kj} - x_{ij}(t)\right)^2 =$$

$$\frac{a(t)^2}{27}\sum_{k=1}^3\left(\mathbb{E}x_{ij}^2(t) - 2p_{kj}\mathbb{E}x_{ij}(t) + \frac{4}{3}p_{kj}^2\right) = \frac{a(t)^2}{9}\left(\mathbb{E}x_{ij}^2(t) - \frac{2}{3}(\sum p_j)\mathbb{E}x_{ij}(t) + \frac{4}{9}\sum p_j^2\right)$$

According to Eq. (4.3), $\mathbb{D}x_{ij}(t) = \mathbb{E}\left(x_{ij}(t) - \frac{1}{3}\sum p_j\right)^2 = \mathbb{E}x_{ij}^2(t) - \frac{2}{3}(\sum p_j)\mathbb{E}x_{ij}(t) +$

$\frac{1}{9}(\sum p_j)^2$.

According to the equations above, we can derive that $\mathbb{D}x_{ij}(t+1) = \frac{a(t)^2}{9}\left(\mathbb{D}x_{ij}(t) + \frac{4}{9}\sum p_j^2 - \frac{1}{9}(\sum p_j)^2\right)$ holds. □

**Proof of corollary 3.1**

The following estimation of $\left|D_{t+1} - \frac{b_{t+1}p_0}{1-b_{t+1}}\right|$ holds:

$\left|D_{t+1} - \frac{b_{t+1}p_0}{1-b_{t+1}}\right| \leq \left|D_{t+1} - \frac{b_t p_0}{1-b_t}\right| + \left|\frac{b_t p_0}{1-b_t} - \frac{b_{t+1}p_0}{1-b_{t+1}}\right| = b_t\left|D_t - \frac{b_t p_0}{1-b_t}\right| + p_0\frac{b_t - b_{t+1}}{(1-b_t)(1-b_{t+1})} \leq \frac{4}{9}\left|D_t - \frac{b_t p_0}{1-b_t}\right| + \frac{36 p_0}{25}\frac{2t+1}{T^2}$

where the second inequality is derived due to the fact that $b_t \leq \frac{4}{9}$. Obviously $\lim_{T \to \infty} \frac{36 p_0}{25}\frac{2t+1}{T^2} = 0$, thus when $T \to \infty$, $\left|D_t - \frac{b_t p_0}{1-b_t}\right| \leq \frac{4}{9}\left|D_{t-1} - \frac{b_{t-1}p_0}{1-b_{t-1}}\right| \leq \cdots \leq \left(\frac{4}{9}\right)^{t-1}\left|D_1 - \frac{b_1 p_0}{1-b_1}\right| \to 0$ holds. □

**Proof of proposition 3.5**

In the GWO, we only need to prove that $\lim_{T \to \infty} \mathbb{D}x_{ij}(T) = 0$. The first-order dynamical system described in Eq. (4.6) can be solved as:

$\mathbb{D}x_{ij}(t) = \mathbb{D}x_{ij}(1)\prod_{\tau=1}^{t-1} b_\tau + p_0 b_{t-1}\left(1 + \sum_{\tau=1}^{t-2}\prod_{\varsigma=\tau}^{t-2} b_\tau\right), t > 1$

Thus,

$\mathbb{D}x_{ij}(T) = \mathbb{D}x_{ij}(1)\prod_{t=1}^{T-1} b_t + p_0 b_{T-1}\left(1 + \sum_{t=1}^{T-2}\prod_{\tau=t}^{T-2} b_\tau\right)$
$\qquad < \mathbb{D}x_{ij}(1) b_{T-1} + p_0 b_{T-1}\left(1 + \sum_{t=1}^{T-2} 1\right)$
$\qquad = \frac{4}{9T^2}\left(\mathbb{D}x_{ij}(1) + p_0(T-1)\right)$

Obviously $\mathbb{D}x_{ij}(1) < +\infty$, thus $\lim_{T \to \infty} \mathbb{D}x_{ij}(T) \leq \lim_{T \to \infty} \frac{4}{9T^2}\left(\mathbb{D}x_{ij}(1) + p_0(T-1)\right) = 0$,

which means $\lim_{t \to \infty} \mathbb{D}x_{ij}(t) = 0$ holds. □